\documentclass[10pt,a4paper]{scrartcl}%

\usepackage[T1]{fontenc}
\usepackage{epsfig}
\usepackage{graphicx}
\usepackage{amssymb}
\usepackage{amsmath}
\usepackage[utf8]{inputenc}
\usepackage{amsfonts}
\usepackage{amsthm}
\usepackage{dsfont}
\usepackage{comment}
\usepackage{enumerate}
\usepackage{mathbbol}
\usepackage{mathrsfs}

\usepackage[a4paper,vmargin={2cm,2cm},hmargin={2.25cm,2.25cm}]{geometry}

\usepackage{hyperref}
\usepackage[english]{babel}

\usepackage{listings}

\usepackage{pgf,tikz,pgfplots}

\pgfplotsset{compat=1.15}
\usetikzlibrary{arrows,shapes,positioning,calc}

\theoremstyle{plain}
\newtheorem{thm}{Theorem}
\newtheorem{lemme}[thm]{Lemma}
\newtheorem{prop}[thm]{Proposition}
\newtheorem{cor}[thm]{Corollary}

\theoremstyle{definition}
\newtheorem{defi}[thm]{Definition}

\theoremstyle{remark}

\newtheorem*{rem}{Remark}

\newcommand{\R}{\ensuremath{\mathbb R}}
\newcommand{\C}{\ensuremath{\mathbb C}}
\newcommand{\D}{\ensuremath{\mathbb D([0,1],\mathbb L)}}
\newcommand{\e}{\ensuremath{\mathbb e}}
\newcommand{\minim}{\ensuremath{\mathcal F}}
\newcommand{\prim}{\ensuremath{\mathcal F^{\downarrow}}}
\newcommand{\antiprim}{\ensuremath{\mathcal F^{\uparrow}}}

\newcommand{\proba}[1]{\ensuremath{\mathbb{P} \left( #1 \right)}}

\newcommand{\luka}{\L ukasiewicz }

\date{2022}

\makeindex

\begin{document}

\title{Random monotone factorisations of the cycle}
\author{Etienne Bellin \footnote{etienne.bellin@polytechnique.edu}
\\
{\normalsize CMAP - Ecole Polytechnique}

}

\maketitle

\begin{abstract}
    In this article we study decreasing and increasing factorisations of the cycle, which are decompositions of the cycle $(1~2\dots n)$ into a product of $n-1$ transpositions satisfying monotonicity conditions. We explicit a bijection between such factorisations and plane trees with $n$ vertices. This will allow us to study some of their combinatorial properties, as well as a geometric representation in terms of laminations, which are non-crossing line segments in the unit disk.
\end{abstract}

\section{Introduction}
For $n \geq 2$ we define the set of \emph{minimal factorisations} of size $n$ by
$$
\minim_n := \left\{ (\tau_1,\dots,\tau_{n-1}) : \forall i ~ \tau_i \text{ is a transposition and } \tau_{n-1} \circ \dots \circ \tau_{1} = (1~2\dots n) \right\}. 
$$

\noindent We also define the set of \emph{decreasing} (resp. \emph{increasing}) factorisations, denoted by $\prim_n$ (resp. $\antiprim_n$), the set of elements in $\minim_n$ with the extra condition that, if $\tau_i = (a_i,b_i)$ with $a_i < b_i$, then the $a_i$'s are in decreasing (resp. increasing) order. For example, $((n ~n-1),(n-1 ~n-2),\dots,(3~2),(2~1))$ is always a decreasing factorisation of $(1 ~2 \dots n)$ ; $((8~9),(8~10),(7~8),(2~3),(2~4),(2~5),(1~2),(1~6),(1~7))$ is a decreasing factorisation of $(1\dots 10)$ and $((1~4),(1~6),(2~4),(3~4),(5~6))$ is an increasing factorisation of $(1\dots 6)$.

\medskip

\noindent Monotone factorisations are part of a wide family of enumerative problems initiated by Hurwitz \cite{Hurwitz}. Indeed the study of factorisations of permutations knows many variants and is linked to many other mathematical fields. For example the number of minimal transitive factorisation is counted by the Hurwitz numbers which is linked to ramified covers of the sphere. We suggest to look at \cite{bousquet} and the references therein for further information about this link and for a generalisation of Hurwitz' theorem to general permutations. Another example is minimal factorisations with transposition of the form $(s,s+1)$, called \emph{reduced decompositions}, which have been studied as well \cite{stanley} \cite{edelman}. These factorisations are also called \emph{sorting networks} since they are linked to sorting algorithms. Factorisations using star transpositions \cite{pak} \cite{irving2009} \cite{goulden2009} \cite{feray} and cycle of given length \cite{biane2004} have been studied as well. Decreasing factorisations are linked to Jucys-Murphy elements and matrix models \cite{matsumoto}. The goal of this article is to study the behaviour of a minimal factorisation chosen uniformly at random in $\prim_n$ or $\antiprim_n$. The asymptotic study of random minimal factorisations has been initiated in \cite{kor2018} and pursued in \cite{kor2019}, \cite{thevenin} and \cite{feray2021}. The study of random sorting networks has also been studied in \cite{angel2007} and \cite{dauvergne}. The present paper completes the previous articles in the case of random decreasing and increasing factorisations.

\medskip

\noindent Decreasing factorisations (also called \emph{primitive factorisations}) have been studied by Gewurz \& Merola \cite{gewurz} who showed that they are counted by Catalan numbers. Two bijections are mentioned in \cite{gewurz}: the first one 
\begin{equation}
\label{eq bij min parking}
((a_1~b_1),\dots,(a_{n-1}~b_{n-1})) \rightarrow (a_{n-1},a_{n-2},\dots,a_1)
\end{equation}
is a bijection with so-called increasing parking functions (i.e.~parking functions $(\pi_1,\dots,\pi_{n-1})$ such that $\pi_1 \leq \dots \leq \pi_{n-1}$), and the second one
\begin{equation}
\label{eq bij min 231}
((a_1~b_1),\dots,(a_{n-1}~b_{n-1})) \rightarrow (b_{n-1},b_{n-2},\dots,b_1)
\end{equation}
is a bijection with 231-avoiding permutations. Increasing factorisations appear in \cite{irving2021} and by the generating function in Theorem 5.2, they are also counted by Catalan numbers.

\medskip

\noindent In this work, we explain how both decreasing and increasing factorisations can be put in bijection with plane trees in a rather unified way. Roughly speaking, it is based on the standard coding of general minimal factorisations by labelled trees. More precisely, a minimal factorisation $((a_1~b_1),\dots,(a_{n-1}~b_{n-1}))$ is encoded by the tree with $n$ vertices, labelled from 1 to $n$, such that an edge labelled $i$ is drawn between the vertices $a_i$ and $b_i$ for all $1 \leq i \leq n-1$. A minimal factorisation can be read on its associated labelled tree, actually, the vertex labels are redundant (see \cite[Lemma~2.6]{kor2019}), meaning that the minimal factorisation can be uniquely retrieved from the edge labels. Our bijections state that the edge labels are also redundant if we are restricted to decreasing or increasing factorisations.

\paragraph{Combinatorial applications.}

Another important merit of our approach is that the bijections allow to obtain information on random uniform decreasing and increasing factorisations by studying the associated random plane tree. For example we show the following results:
\begin{prop}
\label{prop proba prim}
Let $(\tau_1,\dots,\tau_{n-1})$ be a random uniform decreasing factorisation of size $n$ with $\tau_i = (a_i ~ b_i)$ and $a_i < b_i$. Then
\begin{enumerate}
    \item $\{b_1,\dots,b_{n-1}\} = \{2,3,\dots,n\}$ almost surely.
    \item $\proba{\#\{a_1,\dots,a_{n-1}\}=k} = \frac{n+1}{n-1} \binom{n-1}{k} \binom{n-1}{k-1} \binom{2n}{n}^{-1}$ for all $k \in \{1,\dots,n-1\}$.
    \item The convergence 
    $$
    \frac{\#\{a_1,\dots,a_{n-1}\}-n/2}{\sqrt n} \xrightarrow[x\to \infty]{} \mathcal N(0,1/8)
    $$
    holds in distribution where $\mathcal N(0,1/8)$ designates a centered normal distribution with variance $1/8$.
\end{enumerate}
\end{prop}
\begin{rem}
The first item of Proposition \ref{prop proba prim} has already been noticed in \cite{gewurz}. As we said previously, the authors even show that the sequence $(b_1,\dots,b_{n-1})$ is a 231-avoiding permutation, meaning that, if $i<j<k$, then it never happens that $b_k<b_i<b_j$.
\end{rem}
\begin{prop}
\label{prop proba antiprim}
Let $(\tau_1,\dots,\tau_{n-1})$ be a random uniform increasing factorisation of size $n$ with $\tau_i = (a_i ~ b_i)$ and $a_i < b_i$. Then
\begin{enumerate}
    \item The sets $\{a_1,\dots,a_{n-1}\}$ and $\{b_1,\dots,b_{n-1}\}$ are disjoint and their union is $\{1,2,3,\dots,n\}$ almost surely.
    \item Items 2 and 3 of Proposition \ref{prop proba prim} are still satisfied for $\#\{a_1,\dots,a_{n-1}\}$.
\end{enumerate}
\end{prop}
\begin{rem}
Actually, the function given by $(1)$ is also a bijection between increasing factorisations and decreasing parking functions. Moreover one can pass from an increasing to a decreasing parking function just by reversing the order of the elements. Therefore item 2 of Proposition \ref{prop proba antiprim} can be directly deduced from Proposition \ref{prop proba prim}.
\end{rem}

\begin{figure}
\begin{center}
\begin{tikzpicture}
\node[scale =0.4] (a)
{ \scalebox{1}[-1]{\input{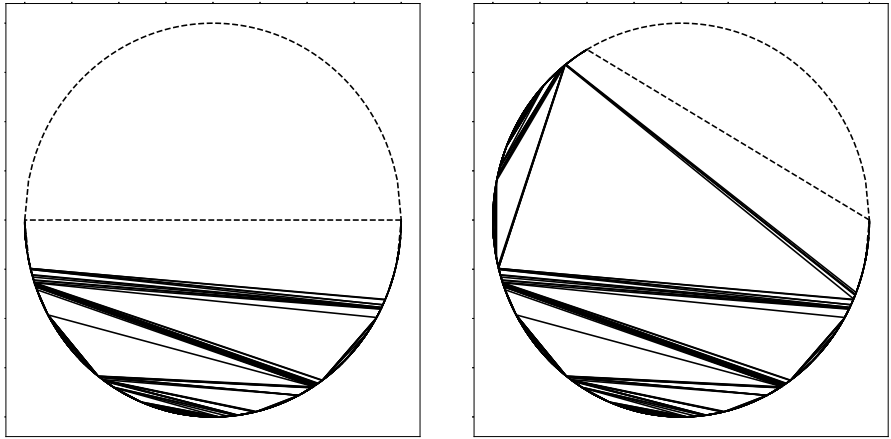}} };
\draw[fill=white,white] (-6.1,-3) rectangle (-5.62,3);
\draw[fill=white,white] (0.1,-3) rectangle (0.53,3);
\draw[fill=white,white] (-6,2.81) rectangle (6,3.1);
\end{tikzpicture}
\caption{A representation of $L_{1/2}(\e)$ on the left and $L_{2/3}(\e)$ on the right.}
\label{fig LBLP}
\end{center}
\end{figure}

\noindent Exploiting bijection (\ref{eq bij min parking}), our work brings some applications to increasing (and decreasing) parking functions. For a parking function $\pi = (\pi_1,\dots,\pi_n)$ of size $n$ and $x \in [0,1]$, set
$$
F_\pi(x) := \frac{1}{n} \#\{i \,:\, \pi_i \leq nx\}.
$$
The following result will be a consequence of Lemma \ref{lemme conv ai} and bijection (\ref{eq bij min parking}).
\begin{prop}
\label{prop fonction de parking}
If $\pi$ is uniformly chosen among the increasing parking functions of size $n$, then the convergence
\begin{equation}
\label{eq conv fonction de parking}
\left(\sqrt{\frac n2}(F_\pi(x)-x)\right)_{0 \leq x \leq 1} \xrightarrow[n\to \infty]{} \e    
\end{equation}
holds in distribution for the uniform norm.
\end{prop}
\begin{rem}
Convergence (\ref{eq conv fonction de parking}), without the $\sqrt 2$, is true when $\pi$ is a uniform minimal factorisation (see \cite[Theorem~15]{diaconis}).
\end{rem}
\paragraph{Geometrical applications.} Our bijections will allow us to study the so-called \emph{lamination process} associated with a random uniform decreasing factorisation. More precisely, let $(\tau_1^n,\dots,\tau^n_{n-1})$ be a random uniform element of $\prim_n$ where $\tau_i^n = (a_i^n~b_i^n)$ with $a_i^n < b_i^n$. The lamination process $(L^n_k)_{1\leq k \leq n-1}$ is defined by
$$
L_k^n := \bigcup_{\ell \leq k} \left[ \exp\left(-2i\pi \frac{a_\ell^n}{n}\right),\exp\left(-2i\pi \frac{b_\ell^n}{n}\right) \right].
$$
It is a process with values in the set of laminations which are unions of non-crossing cords of the unit disk in the complex plane. The main result of this article, besides the bijections, is the following convergence
\begin{thm}
\label{main thm prim}
The convergence 
$$
(L_{\lfloor nt \rfloor}^{n})_{0 \leq t \leq 1} \xrightarrow[n \rightarrow \infty]{} (L_t(\e))_{0 \leq t \leq 1}
$$
holds in distribution in the space $\D$ of càdlàg functions equipped with the Skorokhod's J1 topology.
\end{thm}
\noindent All the topological details and proper definitions can be found in Section \ref{sec asymptotic behaviour} entirely devoted to the proof of Theorem \ref{main thm prim}. The limiting process $L(\e)$ is what we call the \emph{linear Brownian lamination process} and is also properly defined in Section \ref{sec asymptotic behaviour}. The terminology comes from the fact it can be constructed using the normalised Brownian excursion $\e$ (which, informally speaking is a Brownian motion conditioned to reach 0 at time 1 and to be positive on $[0,1]$) by "reading it" from right to left. For $t=1$, $L_1(\e)$ is Aldous' Brownian triangulation \cite{aldous}. Figure \ref{fig LBLP} gives a visual representation of this process. This theorem completes the study of \cite{kor2018} and \cite{thevenin} for random uniform minimal factorisations where the same geometrical point of view has been adopted. Informally speaking, Theorem \ref{main thm prim} encodes how "large" transpositions occur when reading, one after the other, the transpositions of a random uniform decreasing factorisation. To keep this introduction short we postpone discussions about Theorem \ref{main thm prim} to Section \ref{sec discussion} 

\paragraph{Plan of the article.} In Section \ref{sec bij} we construct a bijection between $\prim_n$ (resp. $\antiprim_n$) and the set of plane trees with $n$ vertices. In Section \ref{sec asymptotic behaviour} we study the asymptotic behaviour of a uniform decreasing factorisation of size $n$ when $n$ grows to infinity taking advantage of the bijection found in the previous section in order to prove Theorem \ref{main thm prim}.

\section{Bijections between plane trees, decreasing and increasing factorisations}
\label{sec bij}

\subsection{Labeled trees and plane trees}
\label{sec labeled trees and plane trees}

In this section we recall some basic concepts about trees and introduce some convenient notation that will be used throughout this article.

\medskip

A \emph{tree} is viewed here as a undirected connected acyclic graph. A tree is said to be \emph{labeled} if vertices, or edges, are associated with some real number, called the label of the vertex or the edge. Note that not all the vertices, or edges, of a labeled tree are required to be labeled. We will sometimes emphasize this fact by saying that the tree is \emph{partially labeled}. We will also use the following terminology: a V-labeled tree is a tree where all the vertices, but none of the edges, are labeled, an E-labeled tree is the other way around and an EV-labeled tree is a tree where all the vertices and edges are labeled. If $a$ is a real number and $t$ a labeled tree, we will denote by $v(a,t)$, or simply $v(a)$ if the underlying tree $t$ is unambiguous, the set of vertices of $t$ labeled $a$. If $v(a)$ contains only one vertex, we identify this vertex with the singleton $v(a)$.

\medskip

A \emph{plane tree} is a rooted tree (i.e.~a vertex, called the \emph{root}, is distinguished) where, for every vertex $v$, a total order is given on the edges \emph{stemming from} $v$ (i.e coming from $v$ and going into a child of $v$). Denote by $\Pi_n$ the set of plane trees with $n$ vertices. The labels of a labeled plane tree $t$ are said to be \emph{compatible} with the plane order if for every vertex $v$ of $t$, the natural order of the labels of the edges stemming from $v$ is a restriction of the plane order. Denote by $\Pi_n'$ the set of EV-labeled plane trees with $n$ vertices where the vertices are labeled from $1$ to $n$ (the root being labeled $1$), edges are labeled from $1$ to $n-1$ and such that the labels are compatible with the plane order. In all the figures, when representing a plane tree, the edges stemming from a given vertex will be ordered from left to right. Thus, the smallest edge, for the plane order, will be the first to appear on the left and the biggest will be the furthest on the right. Also, the root will be always represented with a double circle.

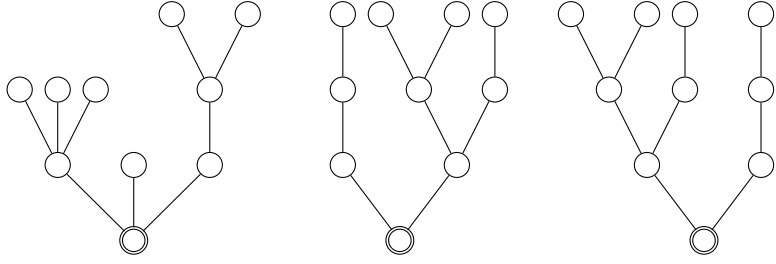
\begin{figure}
\begin{center}
    \begin{tikzpicture}
        \node[draw, double, circle] (1) at (0,0) {};
        \node[draw, circle] (2) at (-1,1) {};
        \node[draw, circle] (3) at (0,1) {};
        \node[draw, circle] (4) at (1,1) {};
        \node[draw, circle] (5) at (-1.5,2) {};
        \node[draw, circle] (6) at (-1,2) {};
        \node[draw, circle] (7) at (-0.5,2) {};
        \node[draw, circle] (8) at (1,2) {};
        \node[draw, circle] (9) at (0.5,3) {};
        \node[draw, circle] (10) at (1.5,3) {};
        \draw (1) -- (2) ;
        \draw (1) -- (3);
        \draw (1) -- (4);
        \draw (2) -- (5);
        \draw (2) -- (6);
        \draw (2) -- (7);
        \draw (4) -- (8);
        \draw (8) -- (9);
        \draw (8) -- (10);
        
        \node[draw, double, circle] (a) at (3.5,0) {};
        \node[draw, circle] (b) at (2.75,1) {};
        \node[draw, circle] (c) at (4.25,1) {};
        \node[draw, circle] (d) at (2.75,2) {};
        \node[draw, circle] (e) at (3.75,2) {};
        \node[draw, circle] (f) at (4.75,2) {};
        \node[draw, circle] (g) at (2.75,3) {};
        \node[draw, circle] (h) at (3.25,3) {};
        \node[draw, circle] (i) at (4.25,3) {};
        \node[draw, circle] (j) at (4.75,3) {};
        \draw (a) -- (b) ;
        \draw (a) -- (c);
        \draw (b) -- (d);
        \draw (c) -- (e);
        \draw (c) -- (f);
        \draw (d) -- (g);
        \draw (e) -- (h);
        \draw (e) -- (i);
        \draw (f) -- (j);
        
        \node[draw, double, circle] (A) at (7.5,0) {};
        \node[draw, circle] (B) at (6.75,1) {};
        \node[draw, circle] (C) at (8.25,1) {};
        \node[draw, circle] (D) at (6.25,2) {};
        \node[draw, circle] (E) at (7.25,2) {};
        \node[draw, circle] (F) at (8.25,2) {};
        \node[draw, circle] (G) at (5.75,3) {};
        \node[draw, circle] (H) at (6.75,3) {};
        \node[draw, circle] (I) at (7.25,3) {};
        \node[draw, circle] (J) at (8.25,3) {};
        \draw (A) -- (B) ;
        \draw (A) -- (C);
        \draw (B) -- (D);
        \draw (B) -- (E);
        \draw (C) -- (F);
        \draw (D) -- (G);
        \draw (D) -- (H);
        \draw (E) -- (I);
        \draw (F) -- (J);
    \end{tikzpicture}
\caption{3 plane trees with 10 vertices. Notice that the last two do not represent the same plane tree. We denote by $\mathfrak{t}$ the tree on the left, it will be our example tree for what follows.} 
\label{fig ex arbre plan et non plan}
\end{center}
\end{figure}

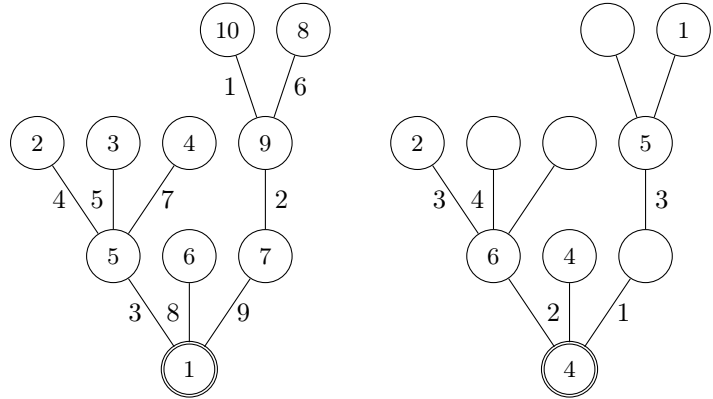
\begin{figure}
\begin{center}
    \begin{tikzpicture}
    [root/.style = {draw,circle,double,minimum size = 20pt,font=\small},
    vertex/.style = {draw,circle,minimum size = 20pt, font=\small}]
        \node[root] (1) at (0,0) {1};
        \node[vertex] (2) at (-1,1.5) {5};
        \node[vertex] (3) at (0,1.5) {6};
        \node[vertex] (4) at (1,1.5) {7};
        \node[vertex] (5) at (-2,3) {2};
        \node[vertex] (6) at (-1,3) {3};
        \node[vertex] (7) at (0,3) {4};
        \node[vertex] (8) at (1,3) {9};
        \node[vertex] (9) at (0.5,4.5) {10};
        \node[vertex] (10) at (1.5,4.5) {8};
        \draw (1) -- node[left] {3} (2) ;
        \draw (1) -- node[left] {8} (3);
        \draw (1) -- node[right] {9} (4);
        \draw (2) -- node[left] {4} (5);
        \draw (2) -- node[left] {5} (6);
        \draw (2) -- node[right] {7} (7);
        \draw (4) -- node[right] {2} (8);
        \draw (8) -- node[left] {1} (9);
        \draw (8) -- node[right] {6} (10);

        \node[root] (a) at (5,0) {4};
        \node[vertex] (b) at (4,1.5) {6};
        \node[vertex] (c) at (5,1.5) {4};
        \node[vertex] (d) at (6,1.5) {};
        \node[vertex] (e) at (3,3) {2};
        \node[vertex] (f) at (4,3) {};
        \node[vertex] (g) at (5,3) {};
        \node[vertex] (h) at (6,3) {5};
        \node[vertex] (i) at (5.5,4.5) {};
        \node[vertex] (j) at (6.5,4.5) {1};
        \draw (a) -- node[left] {} (b) ;
        \draw (a) -- node[left] {2} (c);
        \draw (a) -- node[right] {1} (d);
        \draw (b) -- node[left] {3} (e);
        \draw (b) -- node[left] {4} (f);
        \draw (b) -- node[right] {} (g);
        \draw (d) -- node[right] {3} (h);
        \draw (h) -- node[left] {} (i);
        \draw (h) -- node[right] {} (j);
        
    \end{tikzpicture}
\caption{On the left $\mathfrak{t}$ has been labeled to get $\mathfrak{t}' \in \Pi_{10}'$ (notice that the labels of edges stemming from a same vertex are growing from left to right). On the right $\mathfrak t$ has been partially labeled. This labelling is not compatible with the plane order because the edge label 2 is on the left of 1.} 
\label{fig lexico}
\end{center}
\end{figure}
\medskip

\begin{figure}
\begin{center}
    \begin{tikzpicture}
    [root/.style = {draw,circle,double,minimum size = 20pt,font=\small},
    vertex/.style = {draw,circle,minimum size = 20pt, font=\small}]
        \node[root] (1) at (0,0) {1};
        \node[vertex] (2) at (-1,1.5) {2};
        \node[vertex] (3) at (0,1.5) {6};
        \node[vertex] (4) at (1,1.5) {7};
        \node[vertex] (5) at (-2,3) {3};
        \node[vertex] (6) at (-1,3) {4};
        \node[vertex] (7) at (0,3) {5};
        \node[vertex] (8) at (1,3) {8};
        \node[vertex] (9) at (0.5,4.5) {9};
        \node[vertex] (10) at (1.5,4.5) {10};
        \draw (1) -- (2) ;
        \draw (1) -- (3);
        \draw (1) -- (4);
        \draw (2) -- (5);
        \draw (2) -- (6);
        \draw (2) -- (7);
        \draw (4) -- (8);
        \draw (8) -- (9);
        \draw (8) -- (10);
    \end{tikzpicture}
    \hspace{1em}
    \begin{tikzpicture}[scale=0.5]
        \draw[-stealth] (0,0) -- (11,0);
        \draw[-stealth] (0,-1) -- (0,5);
        
        \draw[thick] (0,0)--(1,2)--(2,4)--(3,3)--(4,2)--(5,1)--(6,0)--(7,0)--(8,1)--(9,0)--(10,-1);
        
        \node[left] (y-1) at (0,-1) {$-1$};
        \node[left] (y0) at (0,0) {$0$};
        \node[left] (y1) at (0,1) {$1$};
        \node[left] (y2) at (0,2) {$2$};
        \node[left] (y3) at (0,3) {$3$};
        \node[left] (y4) at (0,4) {$4$};
        
        \node[below] (x1) at (1,0) {$1$};
        \node[below] (x2) at (2,0) {$2$};
        \node[below] (x3) at (3,0) {$3$};
        \node[below] (x4) at (4,0) {$4$};
        \node[below] (x5) at (5,0) {$5$};
        \node[below] (x6) at (6,0) {$6$};
        \node[below] (x7) at (7,0) {$7$};
        \node[below] (x8) at (8,0) {$8$};
        \node[below] (x9) at (9,0) {$9$};
        \node[above] (x10) at (10,0) {$10$};
        
        \draw[dotted] (x1)--(1,2);
        \draw[dotted] (x2)--(2,4)--(y4);
        \draw[dotted] (x3)--(3,3)--(y3);
        \draw[dotted] (x4)--(4,2)--(y2);
        \draw[dotted] (x5)--(5,1);
        \draw[dotted] (x8)--(8,1)--(y1);
        \draw[dotted] (x10)--(10,-1)--(y-1);
    \end{tikzpicture}
\caption{On the left, the vertices of $\mathfrak{t}$ has been labeled according to their lexicographic order. On the right is represented the \luka path associated with $\mathfrak t$.} 
\label{fig marche de luka}
\end{center}
\end{figure}
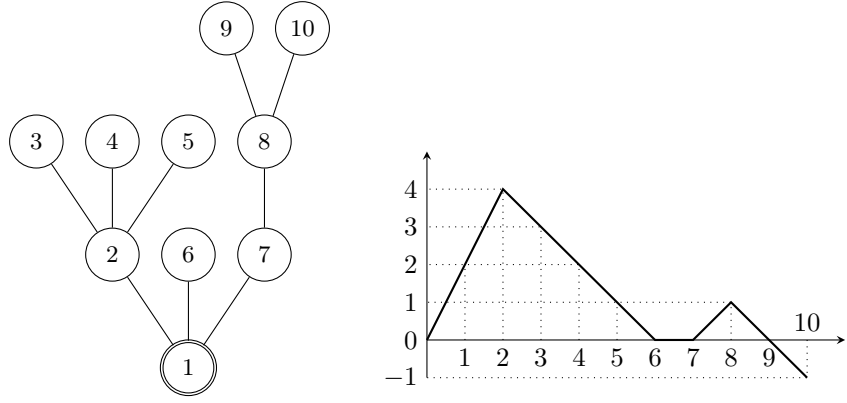

If $t$ is a plane tree then we can define a total order between the vertices of $t$ which is the lexicographic order, also called the depth-first search order, when we encode every vertex $v$ with the list of edges coming from the root and leading to $v$ (see Figure \ref{fig marche de luka} for an example). We denote by $\preceq$ this order and the strict version is denoted by $\prec$. Observe that an ancestor of a vertex $v$ is always smaller than $v$ for the lexicographic order.

\medskip

We end this section with the definition of the \luka path. Let $t$ be a plane tree with $n$ vertices and $v_1,\dots,v_n$ be its vertices ordered in the lexicographic order. The \emph{\luka} path associated with $t$ is the finite sequence $(S_k)_{0\leq k \leq n}$ such that, for all $0 \leq k \leq n$, $S_k+k$ is the number of edges stemming from a vertex $v_\ell$ with $\ell \leq k$ (e.g.~$S_0=0$ and $S_n=-1$). See Figure \ref{fig marche de luka} for a better visualisation.

\subsection{Link between labeled plane trees and minimal factorisations}
\label{sec link plane trees and minim}

\begin{defi}
Let $f$ be a minimal factorisation of size $n$. We define three trees $T_1(f)$, $T_2(f)$ and $T_3(f)$ in the following way:
\begin{itemize}
    \item $T_1(f)$ is the EV-labeled plane tree in $\Pi'_n$ where an edge labeled $k$ is drawn between the vertices labeled $i$ and $j$ iff $\tau_k = (i,j)$.
    \item $T_2(f)$ is the E-labeled plane tree obtained from $T_1(f)$ (we keep the shape, the root and the plane order of $T_1(f)$) by removing the vertex-labels but keeping the edge-labels.
    \item $T_3(f)$ is the (unlabeled) plane tree obtained from $T_1(f)$ by removing the vertex-labels and the edge-labels.
\end{itemize}
\end{defi}
For example, in Figure \ref{fig lexico}, $\mathfrak{t}' = T_1(f)$, where $f=((9~10),(7~9),(1~5),(2~5),(3~5),(8~9),(4~5),(1~6),(1~7))$ is a minimal factorisation of size 10, and $\mathfrak{t}=T_3(f)$. A priori, from the definition above, $T_1(f)$ is only a graph. However it is not hard to see that $T_1(f)$ must be connected, and since it has $n-1$ edges and $n$ vertices, it must be a tree. The map $T_1$ is clearly an injection from the set of all minimal factorisations of size $n$ and the set $\Pi_n'$. However, it is not surjective since the labels of a tree $t \in \Pi'_n$ must satisfy some conditions and can't be chosen freely in order to be in the image of $T_1$. Actually, the latter conditions on the labels are such that, if we erase the vertex-labels of a tree $T_1(f)$, then it is possible to retrieve them. In other words, the map $T_2$ is also an injection. Moreover, any E-labeled plane tree with $n-1$ edges compatibly labeled from 1 to $n-1$ can be uniquely labeled to get an element in the image of $T_1$. We will describe, right after this discussion, an algorithm (see \cite[Section 2.3]{kor2019}) that performs this labeling. To summarize, $T_2$ is a bijection from the set of minimal factorisations of size $n$ to the set of E-labeled plane trees with $n-1$ edges compatibly labeled from 1 to $n-1$. In other words, for minimal factorisations, the vertex-labels are redundant. The goal of the next two sections will be to show that adding the decreasing, or increasing, condition also makes the edge-labels redundant. In other words we will show that $T_3$ is a bijection from the set of decreasing (resp. increasing) factorisations of size $n$ and the set $\Pi_n$ of (unlabeled) plane trees with $n$ vertices.

\medskip 

Let $t$ be a E-labeled plane tree with $n-1$ edges compatibly labeled from 1 to $n-1$. In what follows we describe an "exploration" algorithm \verb&Next&  that takes the tree $t$ and an integer $k \in \{0,1,\dots,n-1\}$ as arguments and gives the label $k+1$ to a vertex of $t$, such that, if we apply successively $\verb&Next&(t,0), \verb&Next&(t,1),\dots,\verb&Next&(t,n-1)$ then $t$ gets the unique labeling that includes him in the image of $T_1$. First $\verb&Next&(t,0)$ gives the label 1 to the root of $t$. For $k\geq 1$, the algorithm starts from the vertex $v_0 = v(k)$ of $t$ labeled $k$. Then it follows the longest possible path possible of edges $e_1,\dots,e_\ell$ such that:
\begin{itemize}
    \item $e_1,\dots,e_\ell$ is a \emph{path of edges} meaning that for all $0<i<\ell$, $e_i$ and $e_{i+1}$ share a common vertex $v_i$.
    \item $e_1$ is the edge with smallest label adjacent to $v_0$.
    \item For all $0<i<\ell$, $e_{i+1}$ is the \emph{successor} of $e_i$ meaning that $e_{i+1}$ has the smallest label among the edges of $v_i$ having a label greater than $e_i$'s label.
\end{itemize}
Starting from $v_0$ and following the path $e_1,\dots,e_\ell$ leads to the vertex of $t$ that will receive the label $k+1$. For instance, if we take the tree $\mathfrak{t}'$ in Figure \ref{fig lexico} again, then $\verb&Next&(\mathfrak{t}',1)$ goes through the edges labeled 3 and 4 to end up on the vertex labeled 2. We can also check that $\verb&Next&(\mathfrak{t}',6)$ passes through the edges labeled 8 and 9 to end up on the vertex labeled 7.

\medskip
 
\noindent Observe that when applying $\verb&Next&(t,0), \verb&Next&(t,1),\dots,\verb&Next&(t,n-1)$ successively, then, after entering the first time in a \emph{fringe subtree} of $t$ (i.e. a subtree of $t$ formed by one vertex of $t$ and all his descendants), the exploration exits this fringe subtree only after labelling all its vertices (see \cite[Proof of Lemma 2.10]{kor2019}). This remark leads to the following fact:

\begin{lemme}
\label{lemme next}
Let $f$ be a minimal factorisation and let $T_1(f)$ be its associated EV-labeled plane tree. For every pair of vertices $u \prec v$ (i.e. $u$ is before $v$ in the lexicographic order) such that $u$ is not an ancestor of $v$, the label of $u$ is smaller than the label of $v$ in $T_1(f)$.
\end{lemme}

Nevertheless, if $u$ is an ancestor of $v$ then the label of $u$ could be greater than the label of $v$. This lemma will be useful in the proof of Proposition \ref{prop bij plane trees and antiprim}.

\subsection{Bijection between plane trees and decreasing factorisations}
\label{sec bij plane trees and prim}

The aim of this section is to show that the map $T_3$, defined in Section \ref{sec link plane trees and minim}, is a bijection between the set $\prim_n$ of decreasing factorisations of size $n$ and the set $\Pi_n$ of plane trees with $n$ vertices. The result is presented in Proposition \ref{prop bij plane trees and prim}, the proof gives an explicit construction of the inverse function of $T_3$. 
\begin{prop}
\label{prop bij plane trees and prim}
The map $T_3$ is a bijection between $\prim_n$ and $\Pi_n$.
\end{prop}

The set of plane trees $\Pi_n$ is known to be of cardinality the $n$-th Catalan's number. So, Proposition \ref{prop bij plane trees and prim}, allows us to deduce the cardinality of $\prim_n$. The cardinality of  $\prim_n$ was already known using bijective (see \cite{gewurz}) and non-bijective (see \cite{irving2021}) proofs, so the result is not new, but it gives a new bijective proof which will be relevant for studying its probabilistic properties.

\begin{cor}
\label{cor card prim}
The cardinality of $\prim_n$ is $\left| \prim_n \right| = \frac{(2n)!}{n!(n+1)!}$.
\end{cor}

\begin{proof}[Proof of Proposition \ref{prop bij plane trees and prim}]
We shall construct explicitly the inverse bijection of $T_3$, i.e.~starting from a rooted plane tree we shall explain how to label its vertices and edges so that it corresponds by $T_3$ to the EV-labeled plane tree of a unique decreasing factorisation. Let $t,t'$ be two partially labeled plane trees, we write $t \sim t'$ if $t$ and $t'$ represent the same plane tree (when we forget all the labels) and for every vertex, or edge, $w$, if $w$ is labeled in both $t$ and $t'$, then it has the same label in $t$ and $t'$. Beware that $\sim$ is not transitive so it's not an equivalence relation. Recall that the map $T_1$ is defined in the beginning of Section \ref{sec link plane trees and minim}. The proof just boils down to showing that for every $t \in \Pi_n$ there exists a unique $f \in \prim_n$ such that $t \sim T_1(f)$.

\medskip

\noindent Let $v_1,\dots,v_n$ be the vertices of $t$ ordered in the lexicographic order (so $v_1$ is the root of $t$) and Let $S$ be the \luka path associated with $t$. Recall that the \luka path is defined in Section \ref{sec labeled trees and plane trees}. For all $k \in \{1,\dots,n\}$ we denote by $t_k$ the partially labeled plane tree obtained from $t$ where:

\begin{itemize}
    \item $\forall i \leq k$, $v_i$ gets the label $i$.
    \item $\forall i \leq k$, the edges stemming from $v_i$ get the labels from $n-1-(S_i+i)$ to $n-1-(S_{i-1}+i-1)$. 
\end{itemize}
In other words, one explores the vertices of $t$ in lexicographic order and labels their outgoing edges in increasing order from left to right using the largest possible labels. We also write $t_0 = t$ (see Figure \ref{fig preuve bij plane trees and prim} for an example).

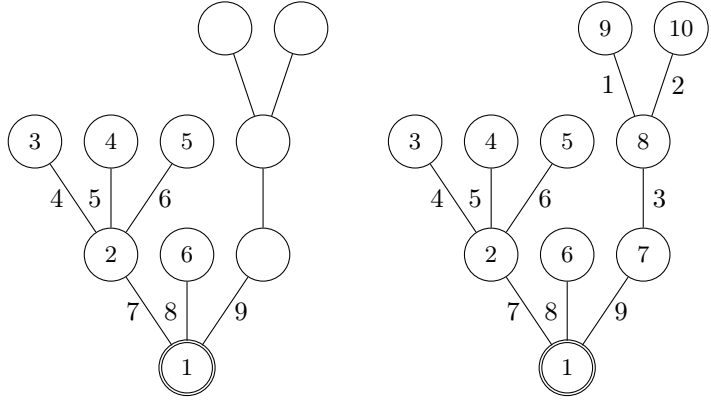
\begin{figure}
\begin{center}
    \begin{tikzpicture}
    [root/.style = {draw,circle,double,minimum size = 20pt,font=\small},
    vertex/.style = {draw,circle,minimum size = 20pt, font=\small}]
        \node[root] (1) at (0,0) {1};
        \node[vertex] (2) at (-1,1.5) {2};
        \node[vertex] (3) at (0,1.5) {6};
        \node[vertex] (4) at (1,1.5) {};
        \node[vertex] (5) at (-2,3) {3};
        \node[vertex] (6) at (-1,3) {4};
        \node[vertex] (7) at (0,3) {5};
        \node[vertex] (8) at (1,3) {};
        \node[vertex] (9) at (0.5,4.5) {};
        \node[vertex] (10) at (1.5,4.5) {};
        \draw (1) -- node[left] {7} (2) ;
        \draw (1) -- node[left] {8} (3);
        \draw (1) -- node[right] {9} (4);
        \draw (2) -- node[left] {4} (5);
        \draw (2) -- node[left] {5} (6);
        \draw (2) -- node[right] {6} (7);
        \draw (4) -- (8);
        \draw (8) -- (9);
        \draw (8) -- (10);

        \node[root] (a) at (5,0) {1};
        \node[vertex] (b) at (4,1.5) {2};
        \node[vertex] (c) at (5,1.5) {6};
        \node[vertex] (d) at (6,1.5) {7};
        \node[vertex] (e) at (3,3) {3};
        \node[vertex] (f) at (4,3) {4};
        \node[vertex] (g) at (5,3) {5};
        \node[vertex] (h) at (6,3) {8};
        \node[vertex] (i) at (5.5,4.5) {9};
        \node[vertex] (j) at (6.5,4.5) {10};
        \draw (a) -- node[left] {7} (b);
        \draw (a) -- node[left] {8} (c);
        \draw (a) -- node[right] {9} (d);
        \draw (b) -- node[left] {4} (e);
        \draw (b) -- node[left] {5} (f);
        \draw (b) -- node[right] {6} (g);
        \draw (d) -- node[right] {3} (h);
        \draw (h) -- node[left] {1} (i);
        \draw (h) -- node[right] {2} (j);
        
    \end{tikzpicture}
\caption{Illustration of the proof of Proposition \ref{prop bij plane trees and prim}. The partially labeled plane tree $\mathfrak{t}_6$ is on the left and the EV-labeled plane tree $\mathfrak{t}_{10}$ is on the right.}
\label{fig preuve bij plane trees and prim}
\end{center}
\end{figure}

\medskip

\noindent Let $f$ be a decreasing factorisation. We will prove by induction that if $f$ is such that $t \sim T_1(f)$, then for all $0\leq k \leq n$, $t_k \sim T_1(f)$.

\medskip

\noindent It is true for $k=0$. Let's assume it is true for a $k \in \{0,\dots,n-1\}$. We advise to look at Figure \ref{fig preuve bij plane trees and prim} while reading the proof to have a better visualisation.

\medskip

\noindent Let $e_1,\dots,e_\ell$ (with $\ell>0$) be the shortest path of edges starting from $v_k$ and going to $v_{k+1}$. Then, let $u_1=v_k,u_2,\dots,u_{\ell+1}=v_{k+1}$ be the successive vertices the path $e_1,\dots,e_\ell$ goes through. Since $v_{k+1}$ is the successor of $v_k$ for the lexicographic order in $t$, one can easily check that: 
\begin{itemize}
    \item $u_1,u_2,\dots,u_\ell$ are the successive ancestors of $v_k$ up to $u_\ell$ in $t$ and $v_{k+1}$ is a child of $u_\ell$.
    \item if $v_k$ is not a leaf, $\ell=1$ and $e_1$ is the smallest edge (for the plane order) among the edges stemming from $v_k$ in $t$.
    \item if $v_k$ is a leaf, $\ell>1$ and $e_\ell$ is the smallest edge among the edges stemming from $u_\ell$ that are greater than $e_{\ell-1}$ in $t$.    
    \item for all $0<i<\ell-1$, $e_i$ is the greatest edge among the edges stemming from $u_{i+1}$ in $t$.
\end{itemize}
Now we observe that if $v$ is a non-root vertex of $t$, then the labels, in $t_n$, of the edges stemming from $v$ are all smaller than the label of the edge between $v$ and its parent. Using this observation, what precedes and by induction we get that $e_1,\dots,e_\ell$ is exactly the path that $\verb&Next&(T_1(f),k)$ would take, so $v_{k+1}$ must have the label $k+1$ in $T_1(f)$. It remains to check that the labels of the edges stemming from $v_{k+1}$ in $T_1(f)$ must be the largest among the unclaimed edge-labels in $t_{k+1}$. Let $e$ be an edge between vertices $u$ and $w$ such that $e$ is not labeled yet in $t_{k+1}$. Denote by $a$ and $b$ the labels of $u$ and $w$ in $T_1(f)$. By construction of $t_{k+1}$, $v_{k+1} \prec u,w$ so $k+1 < a,b$. Since $f$ is a decreasing factorisation, the label of $e$ in $T_1(f)$ must be smaller than the label of the edges stemming from $v_{k+1}$. This implies the initial statement and ends the proof of $t_{k+1} \sim T_1(f)$ which concludes the induction.

\medskip

\noindent Finally we have $T_1(f) = t_n$. By injectivity of $T_1$, uniqueness is shown. It remains to prove existence. Since $T_2$ is surjective, we can find $f$ a minimal factorisation such that $T_2(f)$ is the tree $t_n$ where the vertex-labels are erased. First, notice that, applying $\verb&Next&(T_2(f),0),\dots,\verb&Next&(T_2(f),n-1)$ gives $T_1(f)=t_n$. This can be easily seen if we follow the same scheme as for uniqueness and see that the algorithm $\verb&Next&$ always labels the successor for the lexicographic order in $T_2(f)$. From $T_1(f)$ we easily check that $f$ is a decreasing factorisation.

\end{proof}

The bijection $T_3$ and its inverse function will be useful to prove some asymptotic results in Section \ref{sec asymptotic behaviour} for a decreasing factorisation chosen uniformly at random. The bijection has also direct non-asymptotic consequences presented in Proposition \ref{prop proba prim} which we prove below.

\begin{proof}[Proof of Proposition \ref{prop proba prim}]
Let $t := T_1(\tau_1,\dots,\tau_{n-1})$. For 1, notice that an integer $i$ belongs to $\{b_1,\dots,b_{n-1}\}$ if and only if $v(i,t)$ is adjacent to a vertex whose label is smaller than $i$. Recall that, in $t$, the vertices are labeled according to their lexicographic order. So every vertex $v(i,t)$, for $i>1$, has a parent whose label is smaller than $i$. Consequently $\{b_1,\dots,b_{n-1}\} = \{2,3,\dots,n\}$.

\medskip

\noindent For 2, notice that an integer $i$ belongs to $\{a_1,\dots,a_{n-1}\}$ if and only if $v(i,t)$ is adjacent to a vertex whose label is greater than $i$. This is the case exactly when $v(i,t)$ is not a leaf. Therefore, the cardinality of $\{a_1,\dots,a_{n-1}\}$ is distributed like the number of non-leaf vertices in a uniform plane tree with $n$ vertices. This distribution is known \cite{dershowitz} and is given by Narayana's numbers.

\medskip

\noindent For 3, we apply \cite[Example~2.1]{Janson2016} which gives a central limit theorem for the number of leaves of a critical Bienaymé-Galton-Watson tree with $n$ vertices and reproduction law $\xi$. Taking $\xi$ to a be a geometrical law of parameter $1/2$, gives the law of a uniform plane tree.
\end{proof}

\subsection{Bijection between plane trees and increasing factorisations}
\label{sec bij plane trees and antiprim}

The aim of this section is to show that $T_3$, defined in Section \ref{sec link plane trees and minim}, is a bijection between the set $\antiprim_n$ of increasing factorisations of size $n$ and the set $\Pi_n$ of plane trees with $n$ vertices, following the same principle as in Section \ref{sec bij plane trees and prim}. The result is presented in Proposition \ref{prop bij plane trees and antiprim} and the proof gives an explicit construction of a bijection. 
\begin{prop}
\label{prop bij plane trees and antiprim}
The map $T_3$ is a bijection between $\antiprim_n$ and $\Pi_n$.
\end{prop}

\begin{cor}
\label{cor card antiprim}
The cardinality of $\antiprim_n$ is $\left| \antiprim_n \right|= \frac{(2n)!}{n!(n+1)!}$.
\end{cor}

As for the decreasing case, there are non-bijective proofs of Corollary \ref{cor card antiprim} (see \cite{irving2021}) so this result is not new.

\begin{proof}[Proof of Proposition \ref{prop bij plane trees and antiprim}]

If $t$ is a partially labeled plane tree and $v$ is a vertex of $t$, then we denote by $\verb&sub&(v)$ the fringe subtree of $t$, composed of $v$ and all its descendants, and $\left|\verb&sub&(v)\right|$ the number of vertices of that subtree (e.g. if $v$ is a leaf of $t$, then $\left|\verb&sub&(v)\right|=1$). Similarly to the proof of Proposition \ref{prop bij plane trees and prim}, the proof boils down to show that for every $t \in \Pi_n$ there exists a unique $f \in \antiprim_n$ such that $t \sim T_1(f)$.

\medskip

\noindent To show the latter claim, we will, once again, construct by induction a sequence of partially labeled plane trees $t_1,\dots,t_m$ with $t_0 = t$ and $t_m \in \Pi_n'$ by gradually adding labels on $t$. As before, denote by $v_1,\dots,v_n$ the vertices of $t$ ordered in the lexicographic order. We also denote by $\nu_k$ (resp. $\varepsilon_k$) the smallest unclaimed positive integer vertex (resp. edge) label of $t_k$ (e.g. $\nu_0 = \varepsilon_0 = 1$ and $\nu_m-1=\varepsilon_m = n$). We say that a vertex of $t$ is of \emph{type} 1 if it is a leaf or if its height is an even number (e.g. the root is of type 1 because its height is 0). The other vertices are said to be of \emph{type} 2. To get $t_1$ we simply label the root $v_1$ with $1$ and the edges stemming from $v_1$ with $1,2,\dots,\varepsilon_1-1$ where $\varepsilon_1-1$ is the degree of $v_1$. Now we explain how to construct $t_{k+1}$ from $t_k$. Denote by $v$ the vertex that has been assigned label $\nu_{k-1}$ in $t_k$. If $v = v_n$ then the procedure ends and we set $m = k$. Otherwise, let $u'$ be the successor of $v$ for the lexicographic order and let $u$ be the parent of $u'$. The tree $t_{k+1}$ is constructed in the following way (look at Figures \ref{fig preuve prop bij plane trees and antiprim} and \ref{fig preuve trois cas prop bij plane trees and antiprim}):

\begin{itemize}
\item{Case 1:} If $u'$ is a leaf. We give the label $\nu_k$ to $u'$ and if the edge between $u$ and $u'$ has no label in $t_k$ we label it by $\varepsilon_k$.

\item{Case 2:} If $u'$ is not a leaf and $u$ is of type 2. We give the label $\nu_k$ to $u'$, we give the labels $\varepsilon_k,\dots,\varepsilon_k + j - 1$ to the edges $e_1,\dots,e_j$ stemming from $u'$ and the label $\varepsilon_k + j$ to the edge $e_{j+1}$ between $u'$ and $u$.

\item{Case 3:} If $u'$ is not a leaf and $u$ is of type 1. Denote by $u''$ the successor of $u'$ for the lexicographic order. We give the label $\nu_k + |\verb&sub&(u')| - 1$ to $u'$ and label $\nu_k$ to $u''$. We give the labels $\varepsilon_k,\dots,\varepsilon_k + i - 1$ to the edges $e_1,\dots,e_i$ stemming from $u''$ and the label $\varepsilon_k + i$ to the edge $e_{i+1}$ between $u''$ and $u'$.
\end{itemize}

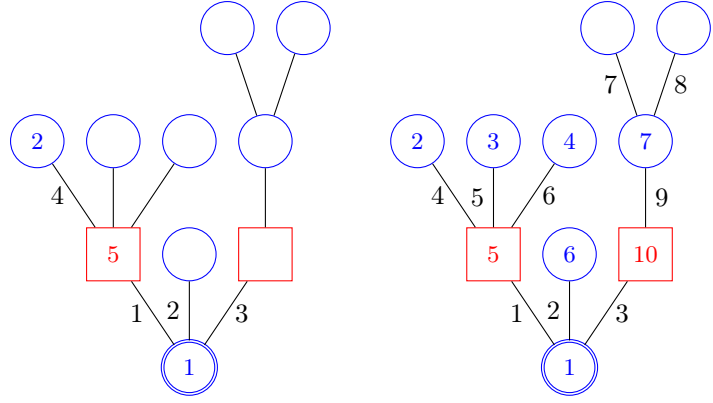
\begin{figure}
\begin{center}
    \begin{tikzpicture}
    [root/.style = {draw,circle,double,minimum size = 20pt,font=\small},
    vertex/.style = {draw,circle,minimum size = 20pt, font=\small}]
        \node[root, blue] (1) at (0,0) {1};
        \node[draw,minimum size = 20pt, font=\small, red] (2) at (-1,1.5) {5};
        \node[vertex, blue] (3) at (0,1.5) {};
        \node[draw,minimum size = 20pt, font=\small, red] (4) at (1,1.5) {};
        \node[vertex, blue] (5) at (-2,3) {2};
        \node[vertex, blue] (6) at (-1,3) {};
        \node[vertex, blue] (7) at (0,3) {};
        \node[vertex, blue] (8) at (1,3) {};
        \node[vertex, blue] (9) at (0.5,4.5) {};
        \node[vertex, blue] (10) at (1.5,4.5) {};
        \draw (1) -- node[left] {1} (2) ;
        \draw (1) -- node[left] {2} (3);
        \draw (1) -- node[right] {3} (4);
        \draw (2) -- node[left] {4} (5);
        \draw (2) -- node[left] {} (6);
        \draw (2) -- node[right] {} (7);
        \draw (4) -- (8);
        \draw (8) -- (9);
        \draw (8) -- (10);

        \node[root, blue] (a) at (5,0) {1};
        \node[draw,minimum size = 20pt, font=\small, red] (b) at (4,1.5) {5};
        \node[vertex, blue] (c) at (5,1.5) {6};
        \node[draw,minimum size = 20pt, font=\small, red] (d) at (6,1.5) {10};
        \node[vertex, blue] (e) at (3,3) {2};
        \node[vertex, blue] (f) at (4,3) {3};
        \node[vertex, blue] (g) at (5,3) {4};
        \node[vertex, blue] (h) at (6,3) {7};
        \node[vertex, blue] (i) at (5.5,4.5) {};
        \node[vertex, blue] (j) at (6.5,4.5) {};
        \draw (a) -- node[left] {1} (b);
        \draw (a) -- node[left] {2} (c);
        \draw (a) -- node[right] {3} (d);
        \draw (b) -- node[left] {4} (e);
        \draw (b) -- node[left] {5} (f);
        \draw (b) -- node[right] {6} (g);
        \draw (d) -- node[right] {9} (h);
        \draw (h) -- node[left] {7} (i);
        \draw (h) -- node[right] {8} (j);
        
    \end{tikzpicture}
\caption{Illustration of the proof of Proposition \ref{prop bij plane trees and antiprim}. The partially lebeled plane tree $\mathfrak{t}_2$ is on the left and $\mathfrak{t}_{6}$ is on the right. Blue circles (resp. red squares) represent vertices of type 1 (resp. 2).} 
\label{fig preuve prop bij plane trees and antiprim}
\end{center}
\end{figure}

We can make useful observations about the construction of the $t_k$'s.
\begin{itemize}
\item{Observation 1:} For every $k$, denote by $i_k$ the number of labeled vertices in $t_k$. The labeled vertices in $t_k$ are exactly the vertices $v_1,\dots,v_{i_k}$. Notice that $i_{k+1} - i_k = 1$ or 2 so in particular $m\leq n$. For every $w, w' \preceq v_{i_k}$ that share a common edge $e$, $e$ is labeled in $t_k$. More precisely, the labeled edges in $t_k$ are exactly the edges adjacent to vertices $w \preceq v_{i_k}$ of type 1.
    
\item{Observation 2:} A vertex $w$ is of type 1 iff there exists $k$ such that $w$ has label $\nu_{k-1}$ in $t_k$. And in this case $w=v_{i_k}$. 
    
\item{Observation 3:} Let $w=v_i$ be a vertex distinct from the root (so $i>1$) or from a leaf with $i\leq i_k$. Denote by $e_1,\dots,e_j$ all the edges stemming from $w$ in increasing order that are labeled in $t_k$ ($j\geq 1$ by observation 1). Also denote by $e$ the edge linking $w$ to its parent (it is labeled in $t_k$ by observation $1$). If $w$ is of type $1$ then $j = \deg(w)-1$ and $e$'s label is greater than $e_j$'s label. If $w$ is of type $2$ then $e$'s label is smaller than $e_1$'s label.

\item{Observation 4:} Let $w,w'$ be two vertices such that $w$ is of type 2 and $w'$ is a descendant of $w$. If $w$ and $w'$ are both labeled in $t_k$ (namely, if $w' \preceq v_{i_k}$) then the label of $w'$ is smaller than the label of $w$.
\end{itemize}

\noindent Now, if $f \in \antiprim_n$ is such that $t \sim T_1(f)$ then we claim that for every $k$, $t_k \sim T_1(f)$. It is clearly true for $k=0$ since $t_0=t$. It is also true for $k=1$ since $f$ is increasing, implying that the transpositions containing 1 must appear first in $f$. Now assume it is true for a certain $k<m$. We use the same notations and the same cases used in the construction of the $t_k$'s. Consider the procedure $\verb&P&$ which consists in applying the algorithm $\verb&Next&$ several times on $T_1(f)$ starting from $v$ (which has label $\nu_{k-1}$ in $T_1(f)$) up to the vertex of label $\nu_k$ in $T_1(f)$. Thus, the procedure $\verb&P&$ consists in applying $\verb&Next&(T_1(f),\nu_{k-1}),
\verb&Next&(T_1(f),\nu_{k-1}+1),\dots,\verb&Next&(T_1(f),\nu_{k}-1)$. We claim that this procedure goes through the vertex $u'$ (but it doesn't mean it stops on $u'$). To see this, we need to study two cases. If $v$ is not a leaf (so $u = v$), since $v$ is of type 1 (observation 2) by using observation 3 we conclude that $\verb&P&$ starts by going up the edge between $v$ and $u'$. If $v$ is a leaf, obviously $\verb&P&$ has no choice but to take the edge between $v$ and its parent. Then observation 3 implies that $\verb&P&$ takes the path down to $u$ (stopping at every vertex of type 2) and then goes to $u'$. The question is now to see if the procedure $\verb&P&$ stops on $u'$ or not.

\begin{itemize}
\item{Case 1:} If $u'$ is a leaf, the procedure stops on $u'$ which imply that $u'$ has label $\nu_k$ in $T_1(f)$. Notice that if $u$ is of type 1, then the edge between $u$ and $u'$ is already labeled in $t_k$ (observation 3) so there is nothing left to prove in order to conclude $t_{k+1} \sim T_1(f)$. Suppose that $u$ is of type 2, we want to show that the edge between $u$ and $u'$ has label $\varepsilon_k$ in $T_1(f)$. Let $e$ be an edge between vertices $v(a)$ and $v(b)$, such that $e$ is not yet labeled in $t_{k+1}$. Suppose that $v(a)\prec v(b)$. Observation 1 implies that, either $u' \prec v(a)$, either $u' \prec v(b)$ and $v(a)$ is an ancestor of $u'$ of type 2. Using observation 4 and Lemma \ref{lemme next}, in both cases $a$ and $b$ are greater than the label of $u'$ in $T_1(f)$, namely $\nu_k$. Since $f$ is an increasing factorisation, it implies that the label of $e$ in $T_1(f)$ must be greater than the label between $u$ and $u'$. It shows that $t_{k+1} \sim T_1(f)$.
    
\item{Case 2:} In this case, $\verb&P&$ stops on $u'$. Indeed, because $u$ is of type 2 it has a greater label than all its descendants (in particular $u'$) in $T_1(f)$ (this comes from the remark preceding Lemma \ref{lemme next} combined with observation 3: applying successively $\verb&Next&$ to $T_2(f)$ will label $u$ after all its descendants). By contradiction, if $\verb&P&$ doesn't stop on $u'$ then, $u'$ has a greater label than its first child $u''$ in $T_1(f)$ (by the remark preceding Lemma \ref{lemme next}). It also implies that the edge between $u$ and $u'$ is smaller than the edge between $u'$ and $u''$. This contradicts the fact  that $f$ is increasing. In conclusion $\verb&P&$ stops on $u'$, thus $u'$ has label $\nu_k$ in $T_1(f)$. It remains to show that the edges adjacent from $u'$ are labeled in $T_1(f)$ as in $t_{k+1}$. Let $e$ be an edge between two vertices labeled $v(a)$ and $v(b)$ such that $e$ is not yet labeled in $t_{k+1}$. Suppose that $v(a) \prec v(b)$. If $u'$ is an ancestor of $v(a)$, since $\verb&P&$ stops on $u'$, $u'$ has a label smaller than $a$ and $b$. If $v(a)$ is an ancestor of $u'$ then $v(a)$ is of type 2 and, like in case 1, we still have that $a,b$ are greater than the label of $u'$. Otherwise, by Lemma \ref{lemme next}, the inequality still holds. Since $f$ is increasing the label of $e$ in $T_1(f)$ must be greater than the labels of the edges adjacent to $u'$ which shows that $t_{k+1}\sim T_1(f)$. 
    
\item{Case 3:} This time $\verb&P&$ doesn't stop at $u'$. We can show this by contradiction, like for case 2. If $\verb&P&$ stops at $u'$ then the edge between $u$ and $u'$ has a greater label than all the edges stemming from $u'$. Moreover it implies $\text{label}(u)<\text{label}(u')<\text{label}(u'')$. Thus the procedure $\verb&P&$ won't stop at $u'$. This shows that $u'$ doesn't have the label $\nu_k$ on $T_1(f)$ and actually has the label $\nu_k + |\verb&sub&(u')| - 1$ (because it is the last vertex seen by $\verb&Next&$ in $\verb&sub&(u')$). From this point we can follow the same reasoning as for case 2 and show that $t_{k+1} \sim T_1(f)$.
\end{itemize} 

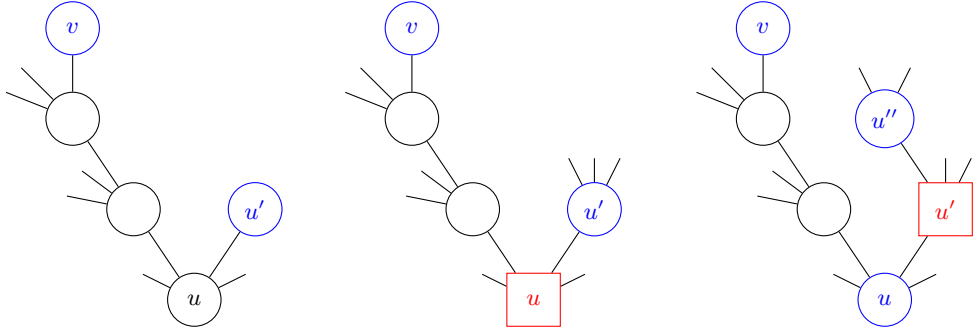
\begin{figure}
\begin{center}    
    \begin{tikzpicture}
    [vertex/.style = {draw,circle,minimum size = 20pt, font=\small},
    scale=0.8]
        \node[vertex] (1) at (0,0) {$u$};
        \node[vertex] (2) at (-1,1.5) {};
        \node[vertex, blue] (3) at (1,1.5) {$u'$};
        \node[vertex] (4) at (-2,3) {};
        \node[vertex, blue] (5) at (-2,4.5) {$v$};
        \node[] (7) at (-3,4) {};
        \node[] (8) at (-3.25,3.5) {};
        \node[] (9) at (-2,2.25) {};
        \node[] (10) at (-2.25,1.75) {};
        \node[] (11) at (-1,0.5) {};
        \node[] (12) at (1,0.5) {};

        \draw (1) -- (2) ;
        \draw (1) -- (3);
        \draw (2) -- (4);
        \draw (4) -- (5);

        \draw (7) -- (4);
        \draw (8) -- (4);
        \draw (9) -- (2);
        \draw (10) -- (2);
        \draw (11) -- (1);
        \draw (12) -- (1);
    \end{tikzpicture}
    \hspace{1 em}
    \begin{tikzpicture}
    [vertex/.style = {draw,circle,minimum size = 20pt, font=\small},
    scale=0.8]
        \node[draw,minimum size = 20pt, font=\small, red] (1) at (0,0) {$u$};
        \node[vertex] (2) at (-1,1.5) {};
        \node[vertex, blue] (3) at (1,1.5) {$u'$};
        \node[vertex] (4) at (-2,3) {};
        \node[vertex, blue] (5) at (-2,4.5) {$v$};
        \node[] (6) at (0.5,2.5) {};
        \node[] (7) at (-3,4) {};
        \node[] (8) at (-3.25,3.5) {};
        \node[] (9) at (-2,2.25) {};
        \node[] (10) at (-2.25,1.75) {};
        \node[] (11) at (-1,0.5) {};
        \node[] (12) at (1,0.5) {}; 
        
        \node[] (13) at (1,2.5) {};
        \node[] (14) at (1.5,2.5) {};

        \draw (1) -- (2) ;
        \draw (1) -- (3);
        \draw (2) -- (4);
        \draw (4) -- (5);
        \draw (3) -- (6);
        \draw (7) -- (4);
        \draw (8) -- (4);
        \draw (9) -- (2);
        \draw (10) -- (2);
        \draw (11) -- (1);
        \draw (12) -- (1);
        \draw (13) -- (3);
        \draw (14) -- (3);

    \end{tikzpicture}
    \hspace{1 em}
    \begin{tikzpicture}
    [vertex/.style = {draw,circle,minimum size = 20pt, font=\small},
    scale=0.8]
        \node[vertex, blue] (1) at (0,0) {$u$};
        \node[vertex] (2) at (-1,1.5) {};
        \node[draw,minimum size = 20pt, font=\small, red] (3) at (1,1.5) {$u'$};
        \node[vertex] (4) at (-2,3) {};
        \node[vertex, blue] (5) at (-2,4.5) {$v$};
        \node[vertex, blue] (6) at (0,3) {$u''$};
        \node[] (7) at (-3,4) {};
        \node[] (8) at (-3.25,3.5) {};
        \node[] (9) at (-2,2.25) {};
        \node[] (10) at (-2.25,1.75) {};
        \node[] (11) at (-1,0.5) {};
        \node[] (12) at (1,0.5) {}; 
        
        \node[] (13) at (1,2.5) {};
        \node[] (14) at (1.5,2.5) {};
        \node[] (15) at (-0.5,4) {};
        \node[] (16) at (0.5,4) {};
        \draw (1) -- (2) ;
        \draw (1) -- (3);
        \draw (2) -- (4);
        \draw (4) -- (5);
        \draw (3) -- (6);
        \draw (7) -- (4);
        \draw (8) -- (4);
        \draw (9) -- (2);
        \draw (10) -- (2);
        \draw (11) -- (1);
        \draw (12) -- (1);
        \draw (13) -- (3);
        \draw (14) -- (3);
        \draw (15) -- (6);
        \draw (16) -- (6);
    \end{tikzpicture}
    
\caption{Illustration of the three cases in the proof of Proposition \ref{prop bij plane trees and antiprim} when $v$ is a leaf. Case 1 is illustrated on the left, case 2 in the middle and case 3 on the right. Black circles represent indifferently vertices of type 1 or 2.} 
\label{fig preuve trois cas prop bij plane trees and antiprim}
\end{center}
\end{figure}

\noindent In conclusion we have $t_k \sim T_1(f)$ for all $k$, so $t_m = T_1(f)$. By injectivity of $T_1$, uniqueness is shown. It remains to prove existence. Since $T_2$ is surjective, we can find $f$ a minimal factorisation such that $T_2(f)$ is the tree $t_m$ where the vertex-labels are erased. First, notice that, applying $\verb&Next&(T_2(f),0),\dots,\verb&Next&(T_2(f),n-1)$ gives $T_1(f)=t_m$. This can be easily seen if we follow the same scheme as for uniqueness and see that the algorithm $\verb&Next&$ labels the vertices of $T_2(f)$ as described above. From $T_1(f)$ we easily check that $f$ is an increasing factorisation.

\end{proof}

Now we can prove the result of Proposition \ref{prop proba antiprim} presented in the introduction.

\begin{proof}[Proof of Proposition \ref{prop proba antiprim}]
Let $t := T_1(\tau_1,\dots,\tau_{n-1})$. Notice that $v(a_1,t),\dots,v(a_{n-1},t)$ are exactly the vertices in $t$ at even height (the root being at height 0). Similarly $v(b_1,t),\dots,v(b_{n-1},t)$ are exactly the vertices in $t$ at odd height. This implies the first point of the corollary. For the second point we use the fact that the number of vertices at even height is distributed like the number of leaves (see \cite{deutsch}). Finally we conclude using the same arguments developed in the proof of Proposition \ref{prop proba prim}.
\end{proof}

\section{Asymptotic behaviour of random uniform decreasing factorisations} \label{sec asymptotic behaviour}

\subsection{Discussing the main result}
\label{sec discussion}

We begin Section \ref{sec asymptotic behaviour} with a discussion on Theorem \ref{main thm prim}.

\medskip

Theorem \ref{main thm prim} is an analog to Theorem 1.2 of \cite{thevenin} where the author studies the lamination process associated to a random uniform minimal factorisation, whereas our theorem concerns random uniform decreasing factorisations. Both theorems give a convergence in the sense of Skorokhod's J1 topology. In the case of minimal factorisations, the convergence of the whole lamination process is a generalisation of the fixed time convergence shown in \cite{kor2018}. In Theorem 1.2 of \cite{thevenin}, the right scaling factor is $\sqrt{n}$ meaning that macroscopic cords appear after a time of order $\sqrt{n}$, whereas in Theorem \ref{main thm prim}, the right scaling factor is $n$. Another difference with the minimal case is that, in the decreasing case, conditionally on \e, the lamination process $L(\e)$ is deterministic while in the minimal case, there is a second layer of randomness. Indeed, conditionally on \e, the cords of the limit process are chosen randomly by throwing points under the curve of $\e$ in a Poissonian way.

\medskip

\noindent Theorem \ref{main thm prim} implies that for any continuous functional $F:\D \rightarrow E$, where $E$ is a metric space, $F(L^n)$ converges towards $F(L(\e))$ in distribution. For instance the functional giving the size of the current longest cord is continuous for the J1 topology (in this case $E$ is the set of càdlàg functions from $[0,1]$ to $\R$). 

\medskip

\noindent We think that there is no result equivalent to Theorem \ref{main thm prim} in the the case of increasing factorisations. More precisely, if $(L_{\lfloor nt \rfloor}^{n})_{0 \leq t \leq 1}$ is the lamination process associated with the random uniform increasing factorisation $(\tau_1^{n},\dots,\tau_{n-1}^{n})$, we believe that there is no convergence in distribution of this lamination process in the sense of Skorokhod's J1 topology. Indeed, using the bijection of Section \ref{sec bij plane trees and antiprim}, one can see that new macroscopic cords could appear at times arbitrarily close in $L^n$, for $n$ large enough, which would make the J1 convergence impossible. More precisely, those macroscopic cords would appear at every branching point of the associated plane tree giving birth to, at least, two subtrees with a non negligible mass of vertices.

\medskip

\noindent To prove Theorem \ref{main thm prim} we introduce, in the next sections, some random excursion $f_n$ associated to the decreasing factorisation $(\tau_1^{n},\dots,\tau_{n-1}^{n})$. There are then two main steps to prove the theorem dealt in the next sections. The first one consists in showing that $L(f_n)$ converges in distribution towards $L(\e)$ for the Skorokhod distance $d_S$. The second step consists in showing that the quasi-distance $d_S'$ (introduced in Section \ref{sec topo}) between $L^n$ and $L(f_n)$ goes to 0 in probability when $n$ tends to infinity. Finally, combining those two steps and using Lemma \ref{lemme almost j1 conv} implies the desired result.

\subsection{Lamination processes}

The aim of this section is to properly introduce the lamination process associated with a random uniform decreasing factorisation of size $n$ (which is quickly defined in the introduction). The ultimate goal being to study the asymptotic behaviour of this process when $n$ tends to infinity. Our main result is Theorem \ref{main thm prim} stating that the lamination process converges in distribution. To make sense of this convergence, first, we need to describe the spaces we are working in and, second, to define the topology we are using. Consider the \emph{unit disk} $D := \{ z\in \C: |z|\leq 1 \}$ of the complex plane.

\begin{defi}[Cord]
Let $u,v \in [0,1]$ with $u\leq v$. The \emph{cord} with end points $u$ and $v$, denoted by $[[u,v]]$, is the subset of the unit disk formed by the segment $\left[e^{-2i\pi u},e^{-2i\pi v} \right]$. More precisely,
$$
[[u,v]] := \left\{ te^{-2i\pi u} +(1-t)e^{-2i\pi v}: t\in[0,1] \right\}.
$$
A cord $[[u,v]]$ is said to be \emph{trivial} if $u=v$ or if $0=u=1-v$. Two cords $[[u,v]]$ and $[[a,b]]$ are said to be \emph{non-crossing} if they do not intersect in the interior of the unit disk, in other words, if $[[u,v]]\cap[[a,b]]\subset \{u,v\}$.
\end{defi}

\begin{defi}[Lamination]
A \emph{lamination} is a non-empty subset $L \subset D$ of the unit disk such that:
\begin{enumerate}[(i)]
    \item $L$ is compact.
    \item $L$ can be written as a union of non-crossing cords.
\end{enumerate}
The set of all laminations is denoted by $\mathbb{L}$.
\end{defi} 

\begin{rem}
The union in (ii) has no particular requirements, it is not necessarily finite or countable. A finite union of non-crossing cords is always a lamination. It is not hard to see that a subset of the unit disk is a lamination if and only if it is the closure of a union of non-crossing cords.
\end{rem}

We can now introduce the object of interest. Let $(\tau_1^{n},\dots,\tau_{n-1}^{n})$ be a random uniform decreasing factorisation of the cycle $(1 \dots n)$. For all $i$, write $(a_i^{n}~b_i^{n})=\tau_i^{n}$ with $a_i^{n} < b_i^{n}$.

\begin{defi}[Discrete lamination process]
\label{def discrete lam proc}
For $1 \leq k \leq n-1$, we define:
$$
L_k^{n} := \bigcup_{i \leq k} \left[\left[\frac{a_i^{n}}{n},\frac{b_i^{n}}{n}\right]\right].
$$
We also set $L_0^{n}:=[[0,0]]=\{(1,0)\}$ and $L_n^{n}:=L_{n-1}^{n}$. The process $L^n := (L^n_{\lfloor nt \rfloor})_{0 \leq t \leq 1}$ is called the \emph{lamination process associated with} $(\tau_1^{n},\dots,\tau_{n-1}^{n})$ or, in short, the \emph{discrete lamination process}, to emphasise the fact that $L^n$ is a finite union of cords appearing at discrete times.
\end{defi}

The denomination "lamination" is justified for $L^n$ since for all $k$, the cords of $L^n_k$ do not cross (see \cite{goulden}). The overall goal of Section \ref{sec asymptotic behaviour} is to show that the process $L^n$ converges in distribution. Below, we describe precisely the limit of $L^n$, which we call the linear Brownian lamination process.

\medskip

\noindent Denote by $\mathcal{E}$ the set of \emph{excursions}, namely, the set of non-negative continuous functions $g$ on $[0,1]$ such that $g(0)=g(1)=0$ and $g(t)>0$ for all $t \in (0,1)$. We say that $[[u,v]]$ is a \emph{cord of} $g \in \mathcal{E}$ if $g(u)=g(v)=\min_{[u,v]}g$. The Brownian excursion, denoted by $\e$, is a random element of $\mathcal{E}$ and will be of particular interest to us. Informally the Brownian excursion is a Brownian motion conditionned to be an element of $\mathcal E$. One can construct $\e$ by renormalising the excursion of a Brownian motion above, or below, 0 around time 1 (see e.g.~\cite{ito}).

\begin{defi}[Linear Brownian lamination process]
\label{def brownian lam proc}
For $g \in \mathcal{E}$ and $t \in [0,1]$ we define:
$$
L_t(g) := \bigcup_{\substack{[[u,v]] \text{ cord of } g \\ u \geq 1-t}} [[u,v]].
$$
We write $L(g) := (L_t(g))_{0 \leq t \leq 1}$. When $g=\e$, $L(\e)$ is called the \emph{linear Brownian lamination process}.
\end{defi}

The limit of $L^n$ is the linear Brownian lamination process $L(\e)$. Proposition \ref{prop L is cadlag} justifies the denomination "lamination" for $L(\e)$ by proving that the cords of $\e$ do not cross and form a compact set. Actually this property is satisfied for any excursion with unique local minima. Formally, an excursion $g$ is said to have \emph{unique local minima} if, whenever $g$ reaches two local minima at distinct times $s$ and $t$, then $g(s) \neq g(t)$.

\subsection{Topology}
\label{sec topo}

In this section we introduce all the topologies and their associated tools to state formally the convergence of Theorem \ref{main thm prim} and prove it. First we need a topology on the set of laminations $\mathbb{L}$. Since laminations are, by definition, non-empty compact subsets of the unit disk $D$, it is natural to endow $\mathbb{L}$ with the Hausdorff distance. Here, we recall the definition of the Hausdorff distance in the general setting.

\begin{defi}[Hausdorff distance]
Let $(E,d)$ be a metric space and denote by $K$ the set of all non-empty, closed and bounded subsets of $E$. The \emph{Hausdorff distance} $d_H$ is a distance on $K$ given by the following: for all $A,B \in K$
$$
d_H(A,B):= \min \{\varepsilon \geq 0: A \subset B^\varepsilon \text{ and } B \subset A^\varepsilon \}
$$
where $U^\varepsilon := \{x\in E \,:\, d(x,U)\leq \varepsilon\}$ is the set of points at distance at most $\varepsilon$ from $U$.
\end{defi}

In our case we consider the Hausdorff distance on the set of non-empty compact subsets of $D$, denoted by $\mathbb K$. In this setting, $(\mathbb K,d_H)$ is a compact metric space and $\mathbb{L}$ is closed so $(\mathbb L,d_H)$ is also compact. Now that $\mathbb L$ is endowed with a metric, we can consider càdlàg functions with values in $\mathbb L$. Recall that $f : [0,1] \rightarrow (E,d)$ is a \emph{càdlàg} function, where $(E,d)$ is a metric space, if for every $t \in [0,1)$, $f$ is continuous on the right of $t$ and for every $t\in(0,1]$, $f$ has a left limit at $t$. The discrete lamination process $L^n$ is obviously càdlàg (its jumps occur at the times $i/n$ for $1\leq i \leq n-1$). The following proposition shows that the linear Brownian lamination process is a càdlàg function with values in $\mathbb L$.

\begin{prop}
\label{prop L is cadlag}
If $g \in \mathcal{E}$ then $L(g) : [0,1] \rightarrow (\mathbb K,d_H)$ is càdlàg on $[0,1]$ and continuous in $1$. Moreover if $g$ has unique local minima then $L_1(g)$ is a lamination (and so is $L_t(g)$ for all $t \in [0,1]$).
\end{prop}

Before proving Proposition \ref{prop L is cadlag} we give a simple notation which will be useful for most of the following proves: we say that two cords $c$ and $c'$ are \emph{$\varepsilon$-close} if $d_H(c,c') \leq 2\pi \varepsilon$. Notice that the cords $[[u,v]]$ and $[[u',v']]$ are $\varepsilon$-close if $|u-u'| \leq \varepsilon$ and $|v-v'| \leq \varepsilon$.

\begin{proof}[Proof of Proposition \ref{prop L is cadlag}]
\begin{itemize}
    \item For $t \in [0,1]$, $L_t(g)$ is compact. Let $(x_n)_n$ be a sequence of $L_t(g)$ that converges to $x$ belonging to the unit disk. For all $n$, denote by $[[u_n,v_n]]$ a cord of $g$ containing $x_n$ with $1-t \leq u_n$. By compacity suppose that $(u_n,v_n)_n$ converges towards $(u,v)$. By continuity $[[u,v]]$ is a cord of $g$ and $1-t \leq u$. Moreover, we can see that $([[u_n,v_n]])_n$ converges to $[[u,v]]$ for the Hausdorff distance which implies that $x$ belongs to the cord $[[u,v]]$.
    
    \item $L(g)$ is right-continuous. Indeed, by contradiction, assume that we can find $\varepsilon >0$ and a sequence $([[u_n,v_n]])_n$ of cords of $g$ such that $(u_n)_n$ increases and converges towards $u$ and for every cord $[[u',v']]$ of $g$ with $u \leq u'$, for all $n$, $[[u_n,v_n]]$ and $[[u',v']]$ are not $\varepsilon$-close. By compacity suppose that $(v_n)_n$ converges towards $v$. Then by continuity, $[[u,v]]$ is a cord of $g$ and $d_H([[u_n,v_n]],[[u,v]]) \rightarrow 0$ when $n \rightarrow \infty$ which yields a contradiction.
    
    \item $L(g)$ has left limits. For all $t \in (0,1]$, we show that the left limit at $t$ is $K_t := \overline{\cup_{s<t}L_s(g)}$. Let $\varepsilon > 0$ and $c_1,\dots,c_n$ be cords of $K_t$ such that for any cord $c$ of $K_t$, there is an $i$ such that $c$ and $c_i$ are $\varepsilon$-close. For all $i$ we can find $s_i < t$ such that there is a cord $c'_i$ of $L_{s_i}(g)$ which is $\varepsilon$-close of $c_i$. Thus for $t>s \geq \max \{s_1,\dots,s_n\}$, $d_H(K_t,L_s(g)) \leq 4 \pi \varepsilon$.
    
    \item $L(g)$ is continuous in 1. By the last point, the left limit at time $1$ is $K_1 = \overline{\cup_{s<1}L_s(g)}$. Since the cord added at time 1 is the trivial cord $[[0,1]]$ which also corresponds to the cord added at time 0, $[[1,1]]$, we have that $\cup_{s<1}L_s(g) = L_1(g)$. By the first point we know that $L_1(g)$ is closed, so $K_1 = L_1(g)$ which shows the continuity at 1.
    
    \item If $g$ has local minima then the cords of $L_1(g)$ do not cross. Indeed if $[[u_1,v_1]]$ and $[[u_2,v_2]]$ are two cords of $g$ with $u_1 \leq u_2$, then it is clear that $u_1 \leq u_2 \leq v_2 \leq v_1$ or $u_1 \leq v_1 \leq u_2 \leq v_2$.
\end{itemize}
\end{proof}

The last task to complete, in order to make the statement of Theorem \ref{main thm prim} fully rigorous, is to define the topology on the set of càdlàg functions from $[0,1]$ to $(\mathbb L,d_H)$ denoted by $\D$. This topology is of course the celebrated Skorokhod's J1 topology (see e.g.~\cite{jacod} chapter VI).

\begin{defi}[Skorokhod distance]
Let $f,g : [0,1] \rightarrow (E,d)$ be two càdlàg functions with values in a metric space $(E,d)$. The \emph{Skorohkod distance} between $f$ and $g$ is
$$
d_S(f,g) := \inf_{\varphi \in \Lambda} \max \left\{ \sup_{[0,1]}|\varphi-Id|\,,\, \sup_{[0,1]} d(f,g\circ \varphi) \right\}
$$
where $\Lambda := \{\varphi:[0,1]\rightarrow[0,1] \,:\, \varphi \text{ is continuous, strictly increasing and bijective}\}$ and $Id \in \Lambda$ is the identity function. 
\end{defi}

The function $d_S$ defined as above is a distance and induces the so called \emph{Skorokhod's J1 topology} on the set of càdlàg functions. This topology is often used to study the convergence of random processes with càdlàg trajectories. At this point, all the tools have been properly defined to fully understand the statement of Theorem \ref{main thm prim}. From here, up to the end of this section, we give some properties of Skorokhod's J1 topology which will be useful to prove Theorem \ref{main thm prim} but are not required to understand the result. Let
$$
\Lambda' := \{\varphi:[0,1]\rightarrow[0,1] \,:\, \varphi \text{ is continuous, non-decreasing and surjective}\}.
$$
Similarly to the Skorokhod distance we define the functions
$$
d_S'(f,g) := \inf_{\varphi \in \Lambda'} \max \left\{ \sup_{[0,1]}|\varphi-Id|\,,\, \sup_{[0,1]} d(f,g\circ \varphi) \right\}
$$
for any càdlàg functions $f,g : [0,1] \rightarrow (E,d)$. Notice that $d_S'$ is not a distance, because it is not symmetric, but it satisfies the triangle inequality as well as the separability axiom. In other words it is a quasi-distance. Obviously we always have
$$
d_S'(f,g) \leq d_S(f,g).
$$
The following lemma shows that the convergence for $d_S'$ is equivalent to the convergence for $d_S$ when $(E,d)$ is compact.

\begin{lemme}
\label{lemme almost j1 conv}
Let $h$ and $(h_n)$ be càdlàg functions on $[0,1]$ with values in a metric space $(E,d)$. Suppose that $(E,d)$ is compact, then
\begin{equation}
\label{eq lemme almost j1 conv}
(d_S'(h,h_n) \xrightarrow[n \rightarrow \infty]{} 0  ~~\text{or}~~ d_S'(h_n,h) \xrightarrow[n \rightarrow \infty]{} 0) ~~\implies~~
d_S(h_n,h) \xrightarrow[n \rightarrow \infty]{} 0.
\end{equation}
\end{lemme}

\begin{proof}
We show that the family $(h_n)$ is relatively compact for the J1 topology using an Arzelà-Ascoli-like criteria for càdlàg functions (see \cite[Theorem 6.3]{ethier}). For a càdlàg function $x$ on $[0,1]$ and $\delta >0$ we define
$$
\omega'(x,\delta) := \inf_{\{t_i\}} \max_i \sup \{d(x(s),x(t)): s,t \in [t_{i-1},t_i) \}.
$$
where the infimum is taken over all the subdivisions $0=t_0 < t_1 < \dots < t_n = 1$ that are $\delta$-\emph{sparse} meaning that for all $i$, $t_i > t_{i-1}+\delta$. By Theorem 6.3 of \cite{ethier} we only need to show that
$$
\lim_{\delta \rightarrow 0} \sup_n \omega'(h_n,\delta) = 0.
$$
Let $\varepsilon > 0$. By Lemma 6.2 of \cite{ethier}, $\lim_{\delta \rightarrow 0} \omega'(h,\delta) = 0$, thus we can find $\delta > 0$ and $\{t_i\}$ a $\delta$-sparse partition such that for all $i$, and $s,t \in [t_{i-1},t_i)$, $d(h(s),h(t)) \leq \varepsilon$. For the rest of the proof, suppose that $\lim_{n\to\infty}d_S(h,h_n)=0$ (the case $\lim_{n\to\infty}d_S(h_n,h)=0$ being similar). Let $(\varphi_n)$ be a sequence of $\Lambda'$ such that
$$
\max \left\{\sup_{[0,1]}|\varphi_n - Id|\,,\, \sup_{[0,1]}d(h,h_n \circ \varphi_n) \right\} \xrightarrow[n \rightarrow \infty]{} 0.
$$
For $n$ large enough, $\{ \varphi_n(t_i) \}$ is $\delta$-sparse and $\sup_{[0,1]}d(h,h_n \circ \varphi_n) \leq \varepsilon$. For $s,t \in [\varphi_n(t_{i-1}),\varphi_n(t_i))$, there exist $s',t' \in [t_{i-1},t_i)$ such that
\begin{align*}
d(h_n(s),h_n(t)) = d(h_n \circ \varphi_n (s'), h_n \circ \varphi_n (t')) & \leq d(h_n \circ \varphi_n (s'), h(s')) + d(h(s'), h(t')) + d(h(t'), h_n \circ \varphi_n (t')) \\
& \leq 3\varepsilon.
\end{align*}
We conclude using the above mentioned criteria of relative compacity. Using (\ref{eq lemme almost j1 conv}) and the relative compacity, we conclude that $h$ is an accumulation point of any sub-sequence of $(h_n)$ for the J1 topology, therefore $(h_n)$ converges towards $h$ for the J1 topology.
\end{proof}

\subsection{A deterministic result}
\label{sec deterministic case}

In this section, we introduce precisely the excursion $f_n$, mentioned in Section \ref{sec discussion}, and show the convergence of $L(f_n)$ towards $L(\e)$. Write $T^n := T_1(\tau_1^{n},\dots,\tau_{n-1}^{n})$ the EV-labeled plane tree associated with the random decreasing factorisation $(\tau_1^{n},\dots,\tau_{n-1}^{n})$. Recall that for all $i$, $v(a_i^{n})$ and $v(b_i^{n})$ are the two vertices adjacent to the edge labeled $i$ in $T^n$. Denote by $(S_k^{n})_{0 \leq k \leq n}$ the \luka walk associated with $T^n$. We complete the trajectories of $S^{n}$ by linear interpolation so that $S^n_t$ is defined for all $t \in[0,n]$. The excursion $f_n$ is then the time renormalization of $S^n$, namely, $f_n(t) := S^n_{(n-1)t}$ for all $t\in[0,1]$. Since $T^n$ is a uniform plane tree with $n$ vertices, $S^n$ has the same law as a random walk, with geometric $\mathcal G(1/2)$ increments, conditioned on reaching $-1$ for the first time at $n$. It is then well known that $f_n$, when correctly renormalised in space, converges towards the Brownian excursion $\e$. More precisely, the convergence
\begin{equation}
\label{eq conv luka to exc bis}
\frac{f_n}{\sqrt{2n}} \xrightarrow[n \rightarrow \infty]{} \e
\end{equation}
holds in distribution for the uniform norm. From the above convergence, one could argue that studying the process $L((2n)^{-1/2}f_n)$ instead of $L(f_n)$ would be more relevant, but one can see that $L(\alpha g)=L(g)$ for every $g\in \mathcal E$ and $\alpha \in \R_+$ so the renormalisation constant can be dropped. Notice that $f_n$ has not necessarily unique local minima, so $L(f_n)$ has no reason to be a lamination-valued process, however it is still a càdlàg function with values in $\mathbb K$. Thus, studying its convergence in the sense of Skorokhod's J1 topology still makes sense.

\medskip

To show that $L(f_n)$ converges in distribution towards $L(\e)$ we will actually show the following deterministic convergence which, roughly speaking, shows a continuity property of the lamination-valued process with respect to the underlying excursion. Combined with (\ref{eq conv luka to exc bis}) and Skorokhod's representation theorem, it implies that $L(f_n)$ converges towards $L(\e)$ in distribution in $\D$. 
\begin{prop}
\label{prop cv deterministe}
Let $g,g_n \in \mathcal E$ such that $(g_n)$ converges to $g$ for the uniform norm and $g$ has unique local minima. Then, the convergence
$$
(L_t(g_n))_{0 \leq t \leq 1} \xrightarrow[n \rightarrow \infty]{} (L_t(g))_{0 \leq t \leq 1}
$$
holds for the Skorokhod's J1 topology.
\end{prop}
 
Before proving Proposition \ref{prop cv deterministe}, we start with a definition and a lemma. 
 
\begin{defi}
\begin{enumerate}[(i)]
\item For $0 \leq u \leq t \leq v \leq 1$, we denote by $[[u,t,v]] := [[u,t]]\cup [[u,v]] \cup [[t,v]]$ the \emph{triangle} formed by the three cords with endpoints $u$, $t$ and $v$.
\item We say that $[[u,t,v]]$ is a triangle of $h \in \mathcal E$ if $[[u,t]]$ and $[[t,v]]$ are both cords of $h$ (in this case $[[u,v]]$ is also a cord of $h$).
\item For $\varepsilon > 0$ we say that the cord $[[u,v]]$ is $\varepsilon$-\emph{big} if $v-u \geq 2\varepsilon$. Similarly, we say that the triangle $[[u,t,v]]$ is $\varepsilon$-big if $[[u,t]]$ and $[[t,v]]$ are both $(\varepsilon/2)$-big.
\item For an $\varepsilon$-big cord $[[u,v]]$ of $h \in \mathcal E$ we define:
$$
\eta_\varepsilon(h,u,v) := \min_{[u+\varepsilon,v-\varepsilon]}h - h(u) \geq 0.
$$
\end{enumerate}
\end{defi}

\begin{lemme}
\label{lemma close cords}
Let $g$ and $(g_n)$ be like in Proposition \ref{prop cv deterministe}. Fix $\varepsilon > 0$. The following assertions are true:
\begin{enumerate}[(i)]
    \item Let $[[u,t,v]]$ be a triangle of $g$ such $u < t < v$. Then, for all $n$, there exists a triangle $[[u_n,t_n,v_n]]$ of $g_n$ such that $\lim_{n\to \infty} (u_n,t_n,v_n) = (u,t,v)$.
    \item Let $\eta >0$, then for $n$ large enough, for every $\varepsilon$-big cord $[[u,v]]$ of $g$ such that $\eta_\varepsilon(g,u,v) \geq \eta$, there exists a cord $[[u_n,v_n]]$ of $g_n$ such that $u_n \in (u,u+\varepsilon]$ and $v_n \in [v-\varepsilon,v)$.
    \item There exists $\eta > 0$ such that for every $\varepsilon$-big cord $[[u,v]]$ of $g$ such that $\eta_\varepsilon(g,u,v) < \eta$, there exists an $\varepsilon$-big cord $[[u',v']]$ of $g$ satisfying $\eta_\varepsilon(g,u',v') = 0$, $u' \in [u,u+\varepsilon]$ and $v' \in [v-\varepsilon,v+\varepsilon]$.
\end{enumerate}
\end{lemme}

\begin{proof}[Proof of Lemma \ref{lemma close cords}]
\begin{enumerate}[(i)]
    \item Let $\delta>0$ such that $u$ and $v$ are not in $[t-\delta,t+\delta]$. Denote by $t_n$ one of the points where $g_n$ reaches its minimum over $[t-\delta,t+\delta]$. Let $u_n := \max\{a \in [0,t-\delta] \,:\, g_n(a) = g_n(t_n)\}$ and $v_n := \min\{b \in [t+\delta,1] \,:\, g_n(b) = g_n(t_n)\}$. Then $[[u_n,t_n,v_n]]$ is a triangle of $g_n$. We can deduce that $(u_n,t_n,v_n)$ converges towards $(u,t,v)$. Indeed, let $(u',t',v')$ be an accumulation point of $(u_n,t_n,v_n)$, then by uniform convergence, $[[u',t',v']]$ is a triangle of $g$ and $t'$ is a minimum of $g$ on $[t-\delta,t+\delta]$. By unicity of local minima $t'=t$. We also have that $u'\leq t-\delta \leq t+\delta \leq v'$. Once again by unicity of local minima, $u'=u$ and $v'=v$.

    \item Fix $n$ such that $||g_n-g||_\infty < \eta/2$. Take $u_n := \max\{a \in [u,u+\varepsilon] \,:\, g_n(a) = g(u)+\eta/2\}$ and $v_n := \min \{b \in [v-\varepsilon,v] \,:\, g_n(b)=g(u)+\eta/2\}$. Then $[[u_n,v_n]]$ is a cord of $g_n$, $u_n \in (u,u+\varepsilon]$ and $v_n \in [v-\varepsilon,v)$.
    
    \item by contradiction, assume that for every $m \geq 1$ we can find an $\varepsilon$-big cord $[[u_m,v_m]]$ of $g$ such that $\eta_\varepsilon(g,u_m,v_m) < 1/m$ and such that for every $\varepsilon$-big cord $[[u',v']]$ of $g$ with $\eta_\varepsilon(g,u',v') = 0$, $u' \notin [u_m,u_m+\varepsilon]$ or $v' \notin [v_m-\varepsilon,v_m+\varepsilon]$. By compacity suppose that $(u_m,v_m)$ converges to $(u',v')$. By uniform convergence $[[u',v']]$ is an $\varepsilon$-big cord of $g$ such that $\eta_\varepsilon(g,u',v')=0$. Thus for $m$ large enough we will have $u' \in [u_m-\varepsilon,u_m+\varepsilon]$ and $v' \in [v_m-\varepsilon,v_m+\varepsilon]$. To avoid contradiction, for $m$ large enough, $u'\in[u_m-\varepsilon,u_m)$. In other words, $(v_m)$ converges towards $v$ and $(u_m)$ converges to the right towards $u'$ without ever touching $u'$. Denote by $t'$ the moment when $g$ reaches its minimum (which is $g(u')$) over $[u'+\varepsilon,v'-\varepsilon]$. For $m$ large enough, $u' < u_m < t' < v_m$. Since $[[u',t']]$ and $[[u_m,v_m]]$ are both cords of $g$, the only possibility is that $u',u_m,t'$ and $v_m$ are minimum of $g$ over $[u',v_m]$ which contradicts the uniqueness of local minima.
\end{enumerate}
\end{proof}

\begin{proof}[Proof of Proposition \ref{prop cv deterministe}]
The idea, roughly speaking, is to perform a time change on $g_n$ that maps all its $\varepsilon$-big triangles to the $\varepsilon$-big triangles of $g$, and to show that the associated lamination-valued processes are close.

\medskip

\noindent The goal is to show that $\lim_{n\to\infty}d_S(L(g_n),L(g))=0$. For every $\varepsilon>0$ and $n$ large enough we will construct $\varphi_n^\varepsilon \in \Lambda$ such that, for all $\varepsilon>0$, $\lim_n||\varphi_n^\varepsilon-Id||_\infty=0$ and for all $\varepsilon>0$ and $n$ large enough, $\sup_{[0,1]}d_H(L_t(g),L_t(g_n\circ \varphi_n^\varepsilon)) \leq 4\pi\varepsilon$. This will show Proposition \ref{prop cv deterministe}. Indeed, notice that for all $t \in [0,1]$, $d_H(L_t(g_n\circ \varphi_n^\varepsilon),L_{1-\varphi_n^\varepsilon(1-t)}(g_n)) \leq 2\pi ||\varphi_n^\varepsilon-Id||_\infty$. Thus $\limsup_n d_S(L(g_n),L(g)) \leq 4\pi\varepsilon$, and since it holds for every $\varepsilon>0$, the convergence follows. 

\medskip

\noindent Now, fix $\varepsilon>0$, we tackle the construction of $\varphi_n^\varepsilon$. Notice that $g$ has a finite number of $\varepsilon$-big triangles denoted by $[[u^1,t^1,v^1]],\dots,[[u^m,t^m,v^m]]$. By Lemma \ref{lemma close cords} (i) we can find triangles $[[u^1_n,t^1_n,v^1_n]],\dots,[[u^m_n,t^m_n,v^m_n]]$ of $g_n$ that converge towards the $\varepsilon$-big triangles of $g$. Let $n_0$ such that for all $n \geq n_0$ and for all $1\leq i \leq m$, $u^i_n \in [0,u^i+\varepsilon]$ and such that $(u^1_n,t^1_n,v^1_n,\dots,u^m_n,t^m_n,v^m_n)$ are ordered in $[0,1]$ like $(u^1,t^1,v^1,\dots,u^m,t^m,v^m)$. To avoid a new notation, just assume that, for all $n \geq n_0$, $u^i_n$ is the maximum element of $[0,u^i+\varepsilon]$ such that $[[u^i_n,t^i_n]]$ is a cord of $g_n$. We define $\varphi_n^\varepsilon$ for all $n \geq n_0$ by setting $\varphi_n^\varepsilon(u^i)=u^i_n$, $\varphi_n^\varepsilon(t^i)=t^i_n$ and $\varphi_n^\varepsilon(v^i)=v^i_n$ for all $1 \leq i \leq m$. We also set $\varphi_n^\varepsilon(0)=0$ and $\varphi_n^\varepsilon(1)=1$ and we complete the definition of $\varphi_n^\varepsilon$ by linear interpolation. Finally, $\varphi_n^\varepsilon$ is an element of $\Lambda$ and clearly $||\varphi_n^\varepsilon-Id||_\infty$ tends to 0. For all $n \geq n_0$, define $g_n^\varepsilon := g_n \circ \varphi_n^\varepsilon$. For all $n \geq n_0$, the function $g_n^\varepsilon$ satisfies the following property: (P1) all $\varepsilon$-big triangle of $g$ is also a triangle of $g_n^\varepsilon$. Let $n_1 \geq n_0$ such that for all $n \geq n_1$, $||\varphi_n^\varepsilon-Id||_\infty \leq \varepsilon/2$. Then for all $n \geq n_1$, $g_n^\varepsilon$ satisfies the following property: (P2) for all $\varepsilon$-big triangle $[[u,t,v]]$ of $g$ there is no element $s \in (u,u+\varepsilon/2]$ such that $[[s,t]]$ is a cord of $g_n^\varepsilon$ (this comes from the maximality of $u_n^i$).

\medskip

\noindent It remains to show that for $n$ large enough $\sup_{[0,1]}d_H(L_t(g),L_t(g_n\circ \varphi_n^\varepsilon)) \leq 4\pi\varepsilon$, equivalently, for $n$ large enough the two following conditions hold: (1) for all cord $[[u,v]]$ of $g$ there is a cord $[[u_n,v_n]]$ of $g_n^\varepsilon$ with $u \geq u_n$ such that $[[u,v]]$ and $[[u_n,v_n]]$ are $(2\varepsilon)$-close and (2) for all cord $[[u_n,v_n]]$ of $g_n^\varepsilon$ there is a cord $[[u,v]]$ of $g$ with $u_n \geq u$ such that $[[u,v]]$ and $[[u_n,v_n]]$ are $(2\varepsilon)$-close. First we focus on proving (1). Let $\eta > 0$ that satisfies (iii) of Lemma \ref{lemma close cords} and $n_2 \geq n_1$ such that, for all $n \geq n_2$, item (ii) of Lemma \ref{lemma close cords} is satisfied with $g_n$ replaced by $g_n^\varepsilon$ (indeed, notice that $g_n^\varepsilon$ converges uniformly towards $g$). Let $n \geq n_2$ and $[[u,v]]$ be a cord of $g$. There are three cases:
\begin{itemize}
    \item If $[[u,v]]$ is not $\varepsilon$-big then $[[u,u]]$ is a cord of $g_n^\varepsilon$ which is $(2\varepsilon)$-close to $[[u,v]]$.
    \item If $[[u,v]]$ is $\varepsilon$-big and $\eta_\varepsilon(g,u,v) \geq \eta$, then by item $(ii)$ of Lemma \ref{lemma close cords} we can find a cord $[[u_n,v_n]]$ of $g_n^\varepsilon$, $\varepsilon$-close from $[[u,v]]$, satisfying $u_n \geq u$.
    \item If $[[u,v]]$ is $\varepsilon$-big and $\eta_\varepsilon(g,u,v) < \eta$, then by item $(iii)$ of Lemma \ref{lemma close cords} we can find an $\varepsilon$-big cord $[[u',v']]$ of $g$, $\varepsilon$-close from $[[u,v]]$, satisfying $u' \geq u$ and $\eta_\varepsilon(g,u',v') = 0$. Thus $[[u',v']]$ is also a cord of $g_n^\varepsilon$ by (P1).
\end{itemize}
Now we focus on (2). By contradiction, assume that (2) fails, so there is an increasing sequence $(n_k)_k$ such that for all $k$, there is a cord $[[u_k,v_k]]$ of $g_{n_k}^\varepsilon$ and no cord $[[u,v]]$ of $g$ such that $u \in [u_k,u_k+2\varepsilon]$ and $v \in [v_k-2\varepsilon,v_k+2\varepsilon]$. In particular for all $k$, $[[u_k,v_k]]$ is $\varepsilon$-big since $[[u_k,u_k]]$ is a cord of $g$. By compacity we can suppose that $(u_k,v_k)$ converges to $(u,v)$. By uniform convergence, $[[u,v]]$ is a cord of $g$ which is $\varepsilon$-big. To avoid any contradiction, $(u_k)$ converges to the right towards $u$. Two cases are possible: $\eta_{\varepsilon}(g,u,v)$ is equal to 0 or not.
\begin{itemize}
    \item If $\eta_{\varepsilon}(g,u,v)>0$, by item $(ii)$ of Lemma \ref{lemma close cords} applied with $g_n$ replaced by $g$ (indeed, the sequence constantly equal to $g$ converges uniformly towards $g$), we can find a cord $[[u',v']]$ of $g$ such that $u'\in (u,u+\varepsilon]$ and $v' \in [v-\varepsilon,v)$. So for $k$ large enough $u' \in [u_k,u_k+2\varepsilon]$ and $v' \in [v_k-2\varepsilon,v_k+2\varepsilon]$ which yields a contradiction.
    \item If $\eta_{\varepsilon}(g,u,v)=0$, let $t \in [u+\varepsilon,v-\varepsilon]$ such that $[[u,t,v]]$ is an $\varepsilon$-big triangle of $g$. By (P1), $[[u,t,v]]$ is also a triangle of $g_{n_k}^\varepsilon$. For $k$ large enough $ u \leq u_k < t < v_k$. Since $[[u_k,v_k]]$, $[[u,t]]$ and $[[t,v]]$ are cords of $g_{n_k}^\varepsilon$, $[[u_k,t,v]]$ is a triangle of $g_n^\varepsilon$. Moreover, for $k$ large enough, $u_k \in [u,u+\varepsilon/2]$ so by (P2), for $k$ large enough $u_k=u$. Finally for $k$ large enough, $u \in [u_k,u_k+2\varepsilon]$ and $v \in [v_k-2\varepsilon,v_k+2\varepsilon]$ which yields a contradiction.
\end{itemize}
\end{proof}

\subsection{A modification of the cord process}

To show the second step of the proof of Theorem \ref{main thm prim}, namely that the quasi-distance $d_S'$ between $L^n$ and $L(f_n)$ goes to 0, we introduce a new cord process $L'(f_n)$ which plays the role of an intermediate between $L^n$ and $L(f_n)$. We start this section with the definition of $L'(f_n)$ and then show that $d_S'(L'(f_n),L(f_n))$ converges towards 0 in probability using, once again, a deterministic result.

\begin{defi}
Let $g \in \mathcal E$ and $[[u,v]]$ a cord of $g$. We say that $[[u,v]]$ is a \emph{maximal cord} of $g$ if there is no cord $[[u',v]]$ of $g$ with $u'<u$ and no cord $[[u,v']]$ of $g$ with $v'>v$. Let $\mathcal M(g)$ be the set of maximal cords of $g$, we define for all $t \in [0,1]$,
$$
L'_t(g) := [[1,1]] \cup \overline{\bigcup_{\substack{[[u,v]] \in \mathcal M(g) \\ u \geq 1-t}} [[u,v]]}
$$
and we write $L'(g) := (L'_t(g))_{0 \leq t \leq 1}$.
\end{defi}

One can show that for any $g \in \mathcal E$ the process $L'(g)$ is càdlàg with values in $\mathbb K$. Indeed one can use the same arguments as in the proof of Proposition \ref{prop L is cadlag} and notice the following fact: if $([[u_n,v_n]])$ is a sequence of maximal cords of $g$ such that $(u_n)$ converges to the right towards $u$ and $(v_n)$ converges towards $v$, then $[[u,v]]$ is either a trivial or a maximal cord of $g$.

\begin{prop}
\label{prop modified cord process}
Let $g,g_n \in \mathcal E$ such that $(g_n)$ converges to $g$ for the uniform norm. Suppose that $g$ has unique local minima and that there is no open interval where $g$ is monotone. Then, the following convergence holds:
$$
d_S'(L'(g_n),L(g_n)) \xrightarrow[n\to\infty]{}0.
$$
\end{prop}

\begin{rem}
The Brownian excursion satisfies almost surely the hypothesis of $g$ in Proposition \ref{prop modified cord process}. This proposition coupled with Skorokhod's representation theorem implies the desired convergence in probability.
\end{rem}

\begin{proof}
Let $\varepsilon >0$ and fix $n$. We will construct a function $\varphi_n^\varepsilon \in \Lambda'$ such that for $n$ large enough, $||\varphi_n^\varepsilon-Id||_\infty \leq \varepsilon$ and $\sup_{t\in[0,1]}d_H(L_t(g_n),L'_{\varphi_n^\varepsilon(t)}(g_n)) \leq 4\pi\varepsilon$. For all $u \in [0,1]$, denote by $\psi_n(u)$ the minimum element of $[0,u]$ such that $[[\psi_n(u),u]]$ is a cord of $g_n$. Then set for all $t \in [0,1]$,
$$
\varphi_n^\varepsilon(t) := 1- (1-t) \wedge \min\{\psi_n(u) \,:\, [[u,s,v]] \text{ is an } \varepsilon\text{-big triangle of }g_n \text{ with }u\geq 1-t\}
$$
with the convention that $\min \emptyset = 1$. Notice that this function is an element of $\Lambda'$. Since $(g_n)$ converges to $g$ and $g$ has unique local minima, for $n$ large enough there is no quadruplet $u,t,v,w$ such that $[[u,t,v]]$ and $[[t,v,w]]$ are $\varepsilon$-big triangles of $g_n$ (this is easily shown by contradiction). From now on suppose that $n$ is large enough so that the previous condition is satisfied. Thus $||\varphi_n^\varepsilon-Id||_\infty \leq \varepsilon$. It remains to prove that for all $t \in [0,1]$, $d_H(L_t(g_n),L'_{\varphi_n^\varepsilon(t)}(g_n)) \leq 4\pi\varepsilon$. Before doing so, notice that, since $g$ is nowhere monotone, for $n$ large enough, there is no interval of length $\varepsilon$ where $g_n$ is monotone. Thus every point of $[0,1]$ is at distance at most $\varepsilon/2$ from a local maximum of $g_n$. Once again, from now one we suppose that $n$ is large enough so that so previous condition is satisfied. 

\medskip

\noindent Let $[[u,v]]$ be a cord of $g_n$. We want to find a maximal cord $[[u',v']]$ of $g_n$ which is $(2\varepsilon)$-close to $[[u,v]]$ such that $u' \geq 1-\varphi_n^\varepsilon(1-u)$. If $[[u,v]]$ is not $\varepsilon$-big and $u\leq 1-\varepsilon$, then there is a local maximum $u' \in [u,u+\varepsilon]$ of $g_n$. The trivial cord $[[u',u']]$ is maximal, $(2\varepsilon)$-close to $[[u,v]]$ and $u' \geq u \geq 1-\varphi_n^\varepsilon(1-u)$. If $[[u,v]]$ is not $\varepsilon$-big and $u> 1-\varepsilon$, then the cord $[[1,1]]$ is $\varepsilon$-close to $[[u,v]]$. Suppose that $[[u,v]]$ is $\varepsilon$-big. If $\eta_\varepsilon(g_n,u,v)>0$, then, adapting Lemma \ref{lemma close cords} (ii), we can find a maximal cord $[[u',v']]$ of $g_n$ which is $\varepsilon$-close to $[[u,v]]$ and such that $u < u' \leq v' < v$. Finally, suppose that $\eta_\varepsilon(g_n,u,v) = 0$. We can find a element $t \in [u+\varepsilon,v-\varepsilon]$ such that $[[u,t,v]]$ is an $\varepsilon$-big triangle of $g_n$. Taking $u' = \psi_n(u)$ and $v'$ the maximum of $[v,1]$ such that $[[u',v]]$ is a cord of $g_n$, we have that $[[u',v']]$ is a maximal cord of $g_n$, $\varepsilon$-close to $[[u,v]]$ and $u' \geq 1-\varphi_n^\varepsilon(1-u)$.

\medskip

\noindent Let $t \in [0,1]$ and $[[u',v']]$ be a maximal cord of $g_n$ such that $u' \geq 1-\varphi_n^\varepsilon(t)$. We want to find a cord $[[u,v]]$, $(2\varepsilon)$-close to $[[u,v]]$ such that $u \geq 1-t$. Suppose that $u' < 1-t$. There is an $\varepsilon$-big triangle $[[u,s,v]]$ of $g_n$ such that $u\geq 1-t$ and $u' \geq \psi_n(u)$. If $g_n(u')>g_n(u)$, then $v' < u$ and the trivial cord $[[u,u]]$ is $\varepsilon$-close to $[[u',v']]$. If $g_n(u')=g_n(u)$, then the cord $[[u,v']]$ is $\varepsilon$-close to $[[u',v']]$.
\end{proof}

\subsection{Link between the cords of the decreasing factorisation and the \luka path}
\label{sec link lam and luka}

In this section we finish the proof of Theorem \ref{main thm prim} by showing that $d_S(L^n,L'(f_n))$ goes to 0 in probability when $n$ tends to infinity. Recall that $(\tau_1^n,\dots,\tau^n_{n-1})$ denotes a uniform decreasing factorisation, that $T^n$ denotes the associated EV-labeled plane tree, $S^n$ its \luka path and $f_n(t) = S^n_{(n-1)t}$.

\begin{prop}
\label{prop derniere etape}
The convergence
$$
d_S(L^n,L'(f_n)) \xrightarrow[n\to\infty]{}0
$$
holds in probability.
\end{prop}

Before proving this result we simplify the process $L'(f_n)$, yet another time, by removing redundant cords of $f_n$. More precisely, a cord $[[u,v]]$ of $f_n$ is said to be $\emph{good}$ if it is a maximal cord such that $v$ is of the form $v=i/(n-1)$ with $i \in \{0,\dots,n-1\}$. Notice that a maximal cord $[[u,v]]$ is good if and only if $f_n(u)$ is an integer. If $[[u,v]]$ is a maximal cord of $f_n$ and $i := \lfloor v(n-1) \rfloor$, then $S^n_i-S^n_{i+1}=-1$, namely the $i$-th step of the \luka path $S^n$ is a downward step. Similarly the last step strictly before $u(n-1)$ is strictly positive. One can easily show that for any maximal cord $[[u,v]]$ of $f_n$ there is a good cord $[[u',v']]$ of $f_n$, $(n-1)^{-1}$-close to $[[u,v]]$ such that $u'\geq u$ (indeed, take $v' = \lfloor v(n-1) \rfloor/(n-1)$ and $u'$ the smallest element such that $[[u',v']]$ is a cord of $f_n$).

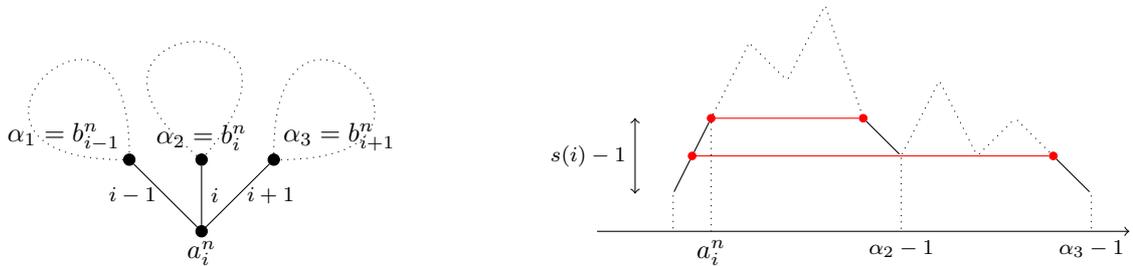
\begin{figure}[!h]
\begin{center}
    \begin{tikzpicture}
    [vertexdot/.style = {draw,circle,fill,inner sep=1.5pt}, scale=0.95]
        \node[vertexdot] (0) at (0,0) {};
        \node[below] at (0,0) {$a^n_i$};
        \node[vertexdot] (1) at (-1,1) {};
        \node[above left] at (-1,1) {$\alpha_1=b^n_{i-1}$};
        \node[vertexdot] (2) at (0,1) {};
        \node[above] at (0,1) {$\alpha_2=b^n_i$};
        \node[vertexdot] (3) at (1,1) {};
        \node[above right] at (1,1) {$\alpha_3=b^n_{i+1}$};
        \draw (0) -- node[left,font=\footnotesize] {$i-1$} (1);
        \draw (0) -- node[right,font=\footnotesize] {$i$} (2);
        \draw (0) -- node[right,font=\footnotesize] {$i+1$} (3);
        \draw[dotted] (1) to[out=180,in=90,distance=3cm] (1);
        \draw[dotted] (2) to[out=140,in=50,distance=3cm] (2);
        \draw[dotted] (3) to[out=90,in=0,distance=3cm] (3);
    \end{tikzpicture}
    \hspace{1em}
    \begin{tikzpicture}[scale=0.5]
    \node[inner sep = 0pt] (0) at (0,0) {};
    \node[draw,circle,fill,inner sep = 1pt,red] (1) at (0.5,1) {};
    \node[draw,circle,fill,inner sep = 1pt,red] (2) at (1,2) {};
    \node[inner sep=0pt] (3) at (2,4) {};
    \node[inner sep=0pt] (4) at (3,3) {};
    \node[inner sep=0pt] (5) at (4,5) {};
    \node[draw,circle,fill,inner sep = 1pt,red] (6) at (5,2) {};
    \node[inner sep = 0pt] (7) at (6,1) {};
    \node[inner sep=0pt] (8) at (7,3) {};
    \node[inner sep=0pt] (9) at (8,1) {};
    \node[inner sep=0pt] (10) at (9,2) {};
    \node[draw,circle,fill,inner sep = 1pt,red] (11) at (10,1) {};
    \node[inner sep = 0pt] (12) at (11,0) {};
    
    \node[below] at (0,-1) {};
    \node[below] at (1,-1) {$a^n_i$};
    \node[below,font=\footnotesize] at (6,-1) {$\alpha_2-1$};
    \node[below,font=\footnotesize] at (11,-1) {$\alpha_3-1$};
    
    \draw (0) -- (1) -- (2);
    \draw[dotted] (2) -- (3) -- (4) -- (5) -- (6);
    \draw (6) -- (7);
    \draw[dotted] (7) -- (8) -- (9) -- (10) -- (11);
    \draw (11) -- (12);
    
    \draw[->] (-2,-1) -- (12,-1);
    
    \draw[<->] (-1,0) -- (-1,2) node[left,pos=0.5,font=\footnotesize] {$s(i)-1$};
    
    \draw[dotted] (0) -- (0,-1);
    \draw[dotted] (2) -- (1,-1);
    \draw[dotted] (7) -- (6,-1);
    \draw[dotted] (12) -- (11,-1);
    
    \draw[red] (2) -- (6);
    \draw[red] (1) -- (11);
    
    \end{tikzpicture}
\caption{Illustration of Proposition \ref{prop derniere etape}. Here $r(i)=2$ and $s(i)=3$. In red are represented the two good cords $c_i'$ and $c'_{i+1}$ of the \luka path $S^n$. Those cords approximate, respectively, the cords $c_i$ and $c_{i+1}$.} 
\label{figure lukasiewicz}
\end{center}
\end{figure}

\begin{figure}[!h]
\begin{center}
    \begin{tikzpicture}
        \node[] (1) at (0,0) {$1$};
        \node[] (2) at (-1,1) {$2$};
        \node[] (4) at (0,1) {$4$};
        \node[] (15) at (1,1) {$15$};
        \node[] (3) at (-2,2) {$3$};
        \node[] (5) at (-1.1,2) {$5$};
        \node[] (10) at (-0.3,2) {$10$};
        \node[] (11) at (0.3,2) {$11$};
        \node[] (12) at (1.1,2) {$12$};
        \node[] (6) at (-2,3) {$6$};
        \node[] (7) at (-1.1,3) {$7$};
        \node[] (8) at (-0.2,3) {$8$};
        \node[] (13) at (0.75,3) {$13$};
        \node[] (14) at (1.5,3) {$14$};
        \node[] (9) at (-0.25,4) {$9$};
        
        \draw (1) --node[left,blue,font=\footnotesize] {12} (2) --node[left,blue,font=\footnotesize] {11} (3);
        \draw (1) --node[left,blue,font=\footnotesize] {13} (4) --node[left,blue,font=\footnotesize] {7} (5) --node[left,blue,font=\footnotesize] {4} (6);
        \draw (5) --node[left,blue,font=\footnotesize] {5} (7);
        \draw (5) --node[right,blue,font=\footnotesize] {6} (8) --node[right,blue,font=\footnotesize] {3} (9);
        \draw (4) --node[left,blue,font=\footnotesize] {8} (10);
        \draw (4) --node[right,blue,font=\footnotesize] {9} (11);
        \draw (4) --node[right,blue,font=\footnotesize] {10} (12) --node[left,blue,font=\footnotesize] {1} (13);
        \draw (12) --node[right,blue,font=\footnotesize] {2} (14);
        \draw (1) --node[right,blue,font=\footnotesize] {14} (15);
        
    \end{tikzpicture}
    \hspace{1em}
    \begin{tikzpicture}[scale = 0.4]
    \draw[->] (0,0) -- (0,7);
    \draw[->] (0,0) -- (16,0);
    
    \node[above,font=\footnotesize] at (0,7) {$S^n$};
    
    \draw (0,0)--(1,2)--(2,2)--(3,1)--(4,4)--(5,6)--(6,5)--(7,4)--(8,4)--(9,3)--(10,2)--(11,1)--(12,2)--(13,1)--(14,0)--(15,-1);
    
    \node[below,font=\footnotesize] at (0,0) {0};
    \node[below,font=\footnotesize] at (1,0) {1};
    \node[below,font=\footnotesize] at (2,0) {2};
    \node[below,font=\footnotesize] at (3,0) {3};
    \node[below,font=\footnotesize] at (4,0) {4};
    \node[below,font=\footnotesize] at (5,0) {5};
    \node[below,font=\footnotesize] at (6,0) {6};
    \node[below,font=\footnotesize] at (7,0) {7};
    \node[below,font=\footnotesize] at (8,0) {8};
    \node[below,font=\footnotesize] at (9,0) {9};
    \node[below,font=\footnotesize] at (10,0) {10};
    \node[below,font=\footnotesize] at (11,0) {11};
    \node[below,font=\footnotesize] at (12,0) {12};
    \node[below,font=\footnotesize] at (13,0) {13};
    \node[below,font=\footnotesize] at (14,0) {14};
    \node[below,font=\footnotesize] at (15,0) {15};
    
    \node[left,font=\footnotesize] at (0,0) {0};
    \node[left,font=\footnotesize] at (0,1) {1};
    \node[left,font=\footnotesize] at (0,2) {2};
    \node[left,font=\footnotesize] at (0,3) {3};
    \node[left,font=\footnotesize] at (0,4) {4};
    \node[left,font=\footnotesize] at (0,5) {5};
    \node[left,font=\footnotesize] at (0,6) {6};
    
    \node[draw,circle,fill,inner sep = 1pt] (1) at (14,0) {};
    \node[draw,circle,fill,inner sep = 1pt] (2) at (12,2) {};
    \node[draw,circle,fill,inner sep = 1pt] (3) at (5,6) {};
    \node[draw,circle,fill,inner sep = 1pt] (4) at (4.5,5) {};
    \node[draw,circle,fill,inner sep = 1pt] (5) at (4,4) {};
    \node[draw,circle,fill,inner sep = 1pt] (6) at (3.67,3) {};
    \node[draw,circle,fill,inner sep = 1pt] (7) at (3.33,2) {};
    \node[draw,circle,fill,inner sep = 1pt] (8) at (1,2) {};
    \node[draw,circle,fill,inner sep = 1pt] (9) at (0.5,1) {};
    \node[draw,circle,fill,inner sep = 1pt] (10) at (0,0) {};
    
    \node[above right,font=\footnotesize] at (1) {$t_{0}$};
    \node[above,font=\footnotesize] at (2) {$t_{2}$};
    \node[above,font=\footnotesize] at (3) {$t_{5}$};
    \node[left,font=\footnotesize] at (4) {$t_{6}$};
    \node[left,font=\footnotesize] at (5) {$t_{8}$};
    \node[right,font=\footnotesize] at (6) {$t_{9}$};
    \node[right,font=\footnotesize] at (7) {$t_{10}$};
    \node[above,font=\footnotesize] at (8) {$t_{13}$};
    \node[right,font=\footnotesize] at (9) {$t_{14}$};
    
    \draw[dotted] (9) -- (13,1) --(13,0);
    \draw[dotted] (7) -- (10,2) -- (10,0);
    \draw[dotted] (6) -- (9,3) -- (9,0);
    \draw[dotted] (4,0)--(5) -- (8,4) -- (8,0);
    \draw[dotted] (4) -- (6,5) -- (6,0);
    
    \end{tikzpicture}
    \hspace{1em}
    \begin{tikzpicture}[scale=3]
    \draw[->] (0,0) -- (0,1.1);
    \draw[->] (0,0) -- (1.1,0);
    
    \node[above,font=\footnotesize] at (0,1.1) {$\varphi_n$};
    \node[right,font=\footnotesize] at (1.1,0) {$t$};
    
    \node[left,font=\footnotesize] at (0,1) {1};
    \node[below,font=\footnotesize] at (1,0) {1};
    \node[below left,font=\footnotesize] at (0,0) {0};
    \node[below,font=\footnotesize] at (10/15,0) {$10/15$};
    \node[left,font=\footnotesize] at (0,10.67/14) {$\frac{16}{21}$};
    
    \draw (0,0)--(2/15,1/7)--(5/15,9/14)--(6/15,9.5/14)--(8/15,10/14)--(10/15,10.67/14)--(13/15,13/14)--(14/15,13.5/14)--(1,1);
    
    \draw[dotted] (0,0)--(1,1);
    \draw[dotted] (0,1)--(1,1)--(1,0);
    \draw[dotted] (10/15,0)--(10/15,10.67/14)--(0,10.67/14);
    \end{tikzpicture}
\caption{Construction of $\varphi_n$ in an example where $n=15$. On the left is represented the random plane tree $T^{15}$. The vertices are labeled in black in the lexicographic order and the edges are labeled in blue. In the middle is drawn the \luka path $S^n$. The abscissa of the black dots are the values of $t_i := (1-\varphi_n(i/n))(n-1)$ and the horizontal dotted lines represent the good cords. On the right is drawn $\varphi_n$.} 
\label{figure ti}
\end{center}
\end{figure}
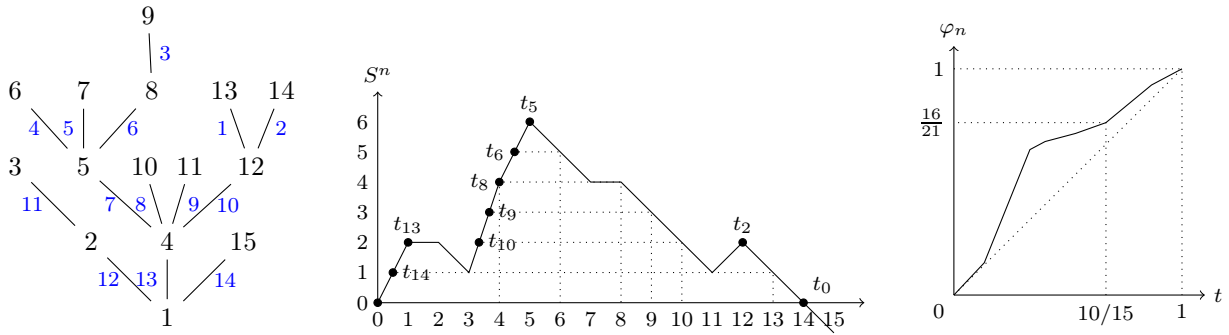

\begin{proof}[Proof of Proposition \ref{prop derniere etape}]
By Skorokhod's representation theorem, suppose that $(f_n)$ converges uniformly towards the Brownian excursion $\e$. With this hypothesis, we will actually show that $d_S(L^n,L'(f_n))$ converges to 0 almost surely. 

\medskip

\noindent For all $1 \leq i \leq n-1$, denote by $c_i$ the $i$-th cord of $L^n$, namely $c_i := [[a^n_i/n,b^n_i/n]]$. For $1 \leq i \leq n-1$, we set $s(i)$ the number of \emph{siblings} of $i$, namely $s(i):= \# \{1 \leq j \leq n-1 \,:\, a^n_j=a^n_i\}$, it is the number of cords having the same left extremity than $c_i$. We also set $r(i)$ the \emph{rank} of $i$, more precisely, $r(i) := \# \{1 \leq j \leq i \,:\, a^n_j=a^n_i\}$, it is the number of cords having the same left extremity than $c_i$ appearing before time $i$. The key to the proof is to notice that for every $i$ of rank $r(i) \geq 2$, the cord $c_i':=[[(a^n_i-h(i))/(n-1),(b^n_i-2)/(n-1)]]$ is a good cord of $f_n$ where $h(i) := (r(i)-2)/(s(i)-1)$. Moreover, all the good cords of $f_n$ are of the form $c'_i$ where $r(i)\geq 2$. This is a consequence of a classical property of the \luka walk (see Figure \ref{figure lukasiewicz}). Indeed, let $\alpha_1,\dots,\alpha_{s(i)}$ be the labels of the children of $a^n_i$ in $T^n$ sorted in increasing order (so in particular $\alpha_1=a^n_i+1$ and $\alpha_{r(i)}=b^n_i$). If we restrict the walk $S^n$ on $[a^n_i,\alpha_{s(i)}-1]$ then the strict descending ladder epochs correspond exactly to the times $a=\alpha_1-1,\dots,\alpha_{s(i)}-1$. Since $h(i)$ is always between 0 and 1, the cord $c'_i$ is $(2(n-1)^{-1})$-close to $c_i$.

\medskip

\noindent As usual we want to find a function $\varphi_n \in \Lambda$ such that $||\varphi_n-Id||_\infty$ and $\sup_{t\in[0,1]} d_H(L^n_{\lfloor nt \rfloor},L'_{\varphi_n(t)}(f_n))$ tend both to 0 as $n$ tends to infinity. What precedes gives a natural choice: for all $1\leq i\leq n-1$ such that $r(i)\geq 2$ we set $\varphi_n(i/n):=1-(a^n_i-h(i))/(n-1)$. We also set $\varphi_n(0):=0$, $\varphi_n(1):=1$ and, as usual, we complete the definition of $\varphi_n$ by linear interpolation (see Figure \ref{figure ti}). By Lemma \ref{lemme conv ai}, who will follow this proof, $a_i^n/n$ gets uniformly close to $1-i/n$, in other words $||\varphi_n-Id||_\infty$ tends to 0 when $n$ goes to infinity. Moreover, it follows from the definition that $\varphi_n$ is strictly increasing.

\medskip

\noindent It remains to show that for every $i \in \{1,\dots,n-1\}$, $d_H(L_i^n,L'_{\varphi^n(i/n)}(f_n))$ tends to 0 as $n$ tends to infinity. Equivalently we need to show that every cord $c_i$ is close to a cord $c'_j$ with $j \leq i$ and every cord $c'_j$ is close to a cord $c_i$ with $i \leq j$. Most of the work has already be done, indeed for all $i$ with $r(i) \geq 2$ the cords $c_i$ and $c'_i$ are $(2(n-1)^{-1})$-close. Hence, it remains only to show that for every cord $c_i$ with $r(i)=1$, there is a close cord $c'_j$ with $j \leq i$. Fix $\varepsilon>0$, as we argued in the proof of Proposition \ref{prop modified cord process}, for $n$ large enough, every point in $[0,1]$ is at distance at most $\varepsilon$ from a local maximum of $f_n$. The final result readily follows, since a local maximum of $f_n$ is a good cord (which is trivial).
\end{proof}

The following lemma shows that $(a^n_{\lfloor nt \rfloor}/n)$ converges uniformly towards $(1-t)$ which is the last ingredient to fully prove Proposition \ref{prop derniere etape}.

\begin{lemme}
\label{lemme conv ai}
The convergence
\begin{equation}
    \label{conv of ai}
    \left( \frac{n(1-t)-a^{n}_{\lfloor nt \rfloor}}{\sqrt{2n}} \right)_{0 \leq t \leq 1} \xrightarrow[n \rightarrow \infty]{} \e
\end{equation}
holds in distribution for the uniform norm where $\e$ is the Brownian excursion on $[0,1]$.
\end{lemme}

\begin{proof}
We use the same notation as in the proof of Proposition \ref{prop derniere etape} and write $a := a^n_i$ to lighten the notations. We show that the following relations hold
\begin{equation}
\label{eq preuve conv of ai}
    S_{a-1}^n \leq n - i - a = S_{a-1}^{n} + (s(i)-r(i)) \leq S_a^{n}.
\end{equation}
For all $1 \leq w \leq n$ denote by $k_w$ the number of children of the vertex labeled $w$ in $T^n$ (so $s(i) = k_a$). The proof of (\ref{eq preuve conv of ai}) is based on the following remark. Let $j\geq i$ then, the edge labeled $j$ stems from a vertex smaller or equal to $a$ for the lexicographic order. Conversely if $w<a$ then all the edges stemming from $w$ have greater labels than $i$. Finally the edges stemming from $a$ that have a greater label than $i$ are exactly those linking $a$ to $\alpha_s$ with $r(i) \leq s \leq s(i)$. Thus we have
    $$
    n-i = \sum_{w<a} k_w + (s(i)-r(i)+1) = S_{a-1}^{n} + a-1 + (s(i)-r(i)+1)
    $$
    and
    $$
    0 \leq s(i)-r(i) \leq s(i)-1.
    $$
This concludes the proof of (\ref{eq preuve conv of ai}). Now we show the convergence (\ref{conv of ai}). Let us rewrite the convergence (\ref{eq conv luka to exc bis}):
\begin{equation}
\left( \frac{S_{nt}^{n}}{\sqrt{2n}} \right)_{0 \leq t \leq 1} \xrightarrow[n \rightarrow \infty]{} \e.
\end{equation}
By Skorokhod's theorem, suppose that the last convergence holds almost surely. Dividing (\ref{eq preuve conv of ai}) by $n$ we deduce that, almost surely 
$$
\left( \frac{a^{n}_{\lfloor nt \rfloor}}{n} \right)_{0 \leq t \leq 1} \xrightarrow[n \rightarrow \infty]{} (1-t)_{0\leq t \leq 1}
$$
for the uniform norm. Dividing (\ref{eq preuve conv of ai}) by $\sqrt{2n}$ then gives (\ref{conv of ai}) (we also use the fact that the Brownian excursion is invariant in law under time inversion).
\end{proof}

\bibliographystyle{alpha}
\bibliography{biblio}

\end{document}